\documentclass[11pt]{article}

\usepackage[utf8]{inputenc}
\usepackage{graphics}
\usepackage{hyperref}
\usepackage{amsmath}
\usepackage{amssymb}
\usepackage{amsfonts}
\usepackage{mathtools}
\usepackage{enumitem}
\usepackage{lipsum}
\usepackage{color,url}
\usepackage{bbm}
\usepackage[normalem]{ulem}
\usepackage{soul}
\usepackage[pagewise]{lineno}

\definecolor{ggreen}{cmyk}{1,     0,      1,      0}
\definecolor{orang}{rgb}{1, 0.40001, 0.2}
\definecolor{viol}{rgb}{0.6, 0.2, 0.9}

\newcommand{\dps}{\displaystyle}
\newcommand{\RR}{\mathbb{R}}

\newcommand{\bfze}{\boldsymbol{0}}
\newcommand{\bfv}{\boldsymbol{v}}
\newcommand{\bfx}{{\boldsymbol{x}}}
\newcommand{\bfu}{\boldsymbol{u}}
\newcommand{\bfy}{\boldsymbol{y}}
\newcommand{\bfz}{\boldsymbol{z}}
\newcommand{\bfJ}{\boldsymbol{J}}

\newcommand{\ba}{\boldsymbol{a}}
\newcommand{\Vv}{\boldsymbol{V}}
\newcommand{\Xx}{\boldsymbol{X}}
\newcommand{\bfB}{\boldsymbol{B}}
\newcommand{\bfC}{\boldsymbol{C}}

\newcommand{\dd}{\mathrm{d}}

\newcommand*{\QEDB}{\hfill\ensuremath{\square}}

\providecommand{\keywords}[1]{\small\textbf{\textit{Keywords:}} #1}
\providecommand{\ams}[1]{\small\textbf{\textit{AMS Subject Classification:}} #1}

\newtheorem{thm}{Theorem}
\newtheorem{prop}[thm]{Proposition}
\newtheorem{lem}[thm]{Lemma}

\newtheorem{rem}[thm]{Remark}

\allowdisplaybreaks

\begin{document}

\baselineskip=15pt


\title{\bf On the Maxwell-Stefan diffusion limit for a reactive mixture of polyatomic gases in non-isothermal setting}

\author{\bf
            B. Anwasia$^{a)}$\footnote{\url{id6226@alunos.uminho.pt}} ,
            M. Bisi$^{b)}$\footnote{\url{marzia.bisi@unipr.it}} ,
            F. Salvarani$^{c)}$\footnote{\url{francesco.salvarani@unipv.it}} ,
            A. J. Soares$^{a)}$\footnote{\url{ajsoares@math.uminho.pt}}
            \\[3mm]
            $^{a)}${\footnotesize
                Universidade do Minho,
                Centro de Matem\'atica,
                Braga, Portugal } \\[1mm]
           $^{b)}${\footnotesize
                Universit\`a di Parma,
                Dipartimento di Scienze Matematiche, Fisiche e Informatiche,
                Parma, Italy } \\[1mm]
           $^{c)}${\footnotesize
                Universit\'e Paris-Dauphine,
                PSL Research University, Ceremade, UMR CNRS
                Paris, France }
                \\
                {\footnotesize \&  Universit\`a degli Studi di Pavia,
                Dipartimento di Matematica,
                Pavia, Italy}}

\date{}

\maketitle


\begin{abstract}
In this article we deduce a mathematical model of Maxwell-Stefan type
for a reactive mixture of polyatomic gases with a continuous structure of internal energy.
The equations of the model are derived in the diffusive limit of a kinetic system of Boltzmann equations
for the considered mixture, in the general non-isothermal setting.
The asymptotic analysis of the kinetic system is performed under a reactive-diffusive scaling
for which mechanical collisions are dominant with respect to chemical reactions.
The resulting system couples the Maxwell-Stefan equations for the diffusive fluxes
with the evolution equations for the number densities of the chemical species
and the evolution equation for the temperature of the mixture.
The production terms due to the chemical reaction and
the Maxwell-Stefan diffusion coefficients are moreover obtained
in terms of the collisional kernels and parameters of the kinetic model.
\end{abstract}

\keywords{
Maxwell-Stefan system,
Reaction-diffusion equations,
Kinetic theory,
Boltzmann equation,
Polyatomic gas mixtures,
Chemical reactions,
Diffusive limit.
}

\ams{
82C40,  
76P05,  
80A32,  
35K57,  
80A30,  
76R50.  
}


\section{Introduction}
\label{sec:int}

Realistic models of multicomponent diffusion phenomena are crucial for many applications in fluid mechanics and chemistry. In the case of isothermal non-reactive gaseous mixtures, composed of at least three different constituents, the diffusive behavior of the species is well described by the equations introduced by Maxwell and Stefan in \cite{max1866, ste1871}, which provide a more general and appropriate framework than the standard Fickian approach \cite{fic1, fic2}.

Despite the popularity of Maxwell-Stefan model and its irreplaceability for many applications in chemical engineering \cite{KW-CES-97, tay-kri-92}, the mathematical properties of the Maxwell-Stefan diffusion equations have only recently been investigated.
In particular, the first studies have been devoted to the matrix formulation of the gradient-flux relationships (see \cite{gio_book} and the references therein).
Subsequently, the well-posedness and the long-time behavior of the solutions of the Maxwell-Stefan system (or some variants) have been studied in \cite{bot-11, bou-gre-sal-12, che-jun-15, HERBERG2017264, hut-sal-18, jun-ste-13}, the numerical simulation of the Maxwell-Stefan system has been the subject of \cite{bon-bou-gre-18, bou-gre-sal-12, gei-15, mcl-bou-14}, and the relationships between Fickian diffusion and the Maxwell-Stefan model have been analyzed in \cite{BGP-NA-17, sal-soa-18}.
By following the research line initiated by Bardos, Golse and Levermore in the Nineties -- whose goal was the derivation of the equations of fluid mechanics starting from the Boltzmann equation \cite{BGL91, BGL93} -- {several} articles have carried out the formal derivation of isothermal multicomponent Maxwell-Stefan type diffusion equations starting from the Boltzmann system for monatomic non-reactive gaseous mixtures \cite{BGP-NA-17, bou-gre-pav-sal-14, bou-gre-sal-15, HS-NA-17, HS-MMAS-17}.
The diffusive limit {in a reactive mixture described by} the simple reacting sphere kinetic model (SRS), which retains the main features of the reaction mechanism without taking into account the internal degrees of freedom of the particles, has been investigated in \cite{AGS-17}.

However, if we consider gaseous mixtures composed of polyatomic gases, with vibrational and rotational degrees of freedom, the standard monatomic models are {no} longer valid and the presence of possible chemical reactions in the mixture can considerably modify the behavior of the system, making the model non-isothermal and non-conservative (at the level of the single species).
The derivation of the precise structure of the equations describing the diffusive regime in this situation is -- of course -- crucial, and this is precisely the goal of the present article.

Our approach consists in obtaining the macroscopic equations starting from a kinetic system of equations defined in the phase space, under the diffusive scaling. We treat the effects of the binary interactions between particles as a simple scattering event involving, at the microscopic scale, only some fundamental laws (in particular, the conservation of momentum and of the total energy), without needing supplementary phenomenological hypotheses.
Because of the solidity of the kinetic structure, the main macroscopic collective features of the system can subsequently be rigorously deduced by means of an appropriate limiting procedure of the kinetic model, rather than heuristically introduced in the macroscopic model.
The result of this approach is a hierarchy of simpler models, which are validated by a precise asymptotic study considering the order of magnitude of the relevant parameters, and are suitable to be applied {in} some particular regime.

In the literature, there is a {variety of kinetic models} which have been proposed to go beyond the monatomic and non-reactive setting.
Polyatomic non reactive mixtures have been studied at the kinetic level, for example,
in \cite{bis-rug-spi-18, bob-bis-gro-spi-18, bou-des-let-per-94, des-97},
whereas various kinetic models for polyatomic reactive mixtures have been derived in
{\cite{bis-cac-16, bis-mon-soa-18, DMS-EJM-05, gro-spi-99}.}

The polyatomic structure of particles may be modeled by means of a set of discrete internal energy levels, or through a proper continuous internal energy variable. The basic features of the discrete levels description may be found in~\cite{gio_book}, while the state of art of kinetic and Extended {Thermodynamics} approaches with a continuous energy have been summarized in~\cite{Ruggeri-Sugiyama}.

In this article we consider, as starting point, the kinetic system proposed in \cite{DMS-EJM-05}, based on the Borgnakke-Larsen procedure \cite{bor-lar-75}, because of its main features and advantages detailed in~\cite{DMS-EJM-05}.
The model introduced in \cite{DMS-EJM-05} describes indeed a mixture of reactive polyatomic gases by adding to the usual independent variables of the phase-space of the system (time $t$, position $\bfx$ and velocity $\bfv$) a continuous internal energy variable $I\in\mathbb{R}^+$, which governs, together with the kinetic energy, the binary encounters -- both of reactive and of non-reactive type.
By carefully choosing a set of measures $\varphi_i(I)\mathrm{d} I$, the model is moreover consistent, at the macroscopic level, with the energy law of any type of polyatomic {gas}, and it does not need to take into account a large number of discrete energy levels.

In our article, for the sake of clarity, we limit our study to a quaternary mixture of polyatomic gases in the presence of
a reversible chemical reaction of bimolecular type. This framework allows us to work in the general non-isothermal setting and to derive from the kinetic equations, at the formal level,
a coupled system of equations which governs the diffusion phenomena in the mixture in the presence of chemical reactions.

More specifically, the set of {equations obtained in this way} include the evolution equations for the chemical species,
the Maxwell-Stefan equations for the diffusive fluxes and the evolution equation for the temperature of the mixture.
Of course, generalizations to more complicated mixtures are possible by introducing straightforward modifications in the computations.

With respect to the standard isothermal non-reactive Maxwell-Stefan system, the equations obtained in this article show some crucial differences.
First of all, the continuity equations for the various species, which would assure the conservation of the molar densities of the species, are replaced by balance equations, whose right-hand side takes into account the effects of the chemical reactions on the densities of the reactants.
The balance terms, {once the equilibrium is reached, guarantee the validity of the law of mass action,} which depends on the internal energy structure of the species, on the temperature of the mixture and on the {reaction heat $E$.}
We need moreover to take into account the energy balance due to the effects of the chemical reactions.

We highlight that the target equations derived in this article are, up to the best of our knowledge, the only ones which take into account both the polyatomic structure of the constituents of the mixtures and the effect of chemical reactions in the non-isothermal setting.

The structure of the article is the following. After introducing, in Section \ref{sec:meqs}, the model governing the reactive mixture of polyatomic gases proposed in \cite{DMS-EJM-05} and its main properties (conservation laws, equilibrium states, chemical rates and H-theorem), we consider in Section \ref{sec:sea} the scaled
system in the diffusive regime.
Finally, the limiting diffusive equations for the reactive mixture are deduced in Section~\ref{sec:limin}.
The conclusion, in Section \ref{sec:concl}, which summarizes our results, is followed by an appendix which gathers some technical computations that are necessary for deducing the results of Section \ref{sec:limin}.


\section{The model for a reactive mixture of polyatomic gases}
\label{sec:meqs}
In this section we briefly present the kinetic model proposed {in \cite{DMS-EJM-05}}
for a quaternary reactive mixture of polyatomic gases
with a continuous structure of internal energy.
We restrict our presentation to the kinetic equations, the central aspects of the collisional dynamics,
conservation laws, equilibrium states and $\cal H$-theorem.

\medskip


\noindent
Following \cite{DMS-EJM-05}, we consider a quaternary mixture of species $A_1,A_2,A_3$ and $A_4$,
participating in the reversible chemical reaction of type
\begin{equation}
A_1 + A_2 \leftrightharpoons A_3 + A_4 \, .
\label{eq:cr}
\end{equation}
For each species $A_i$, {with $i=1,2,3,4,$} we introduce its distribution function $f_i$
which depends on time $t\in\RR^+$, position $\bfx \in \RR^3$, velocity $\bfv \in \RR^3$
and on the internal energy variable $I \in [0,+\infty[$.
For sake of simplicity, {in many cases we omit} the dependence of each $f_i$ on $t$ and $\bfx$,
and write $f_i(\bfv,I)$. Sometimes we simply write $f_i$.

We denote the molecular mass of each species by $m_i$ and the chemical binding energy by $E_i$.
Furthermore, we introduce a weight $\varphi_i(I)$, which aims at obtaining the energy law of polyatomic gases and
the mass action law of chemical kinetics.
With reference to the chemical reaction \eqref{eq:cr},
the conservation of mass requires that
\begin{equation}
m_1 + m_2 = m_3 + m_4 = M ,
\label{eq:conservation_mass_chem}
\end{equation}
and the balance of binding energies is specified by the reaction heat
\begin{equation}
E = E_3 + E_4 - E_1 - E_2 ,
\label{eq:E}
\end{equation}
such that $E > 0$ means that the forward reaction $A_1 + A_2 \rightarrow A_3 + A_4$ is endothermic
whereas $E < 0$ indicates that it is exothermic.

\medskip


\noindent
An important aspect of the present description is that the moments of the distribution function $f_i(t,\bfx,\bfv,I)$
are defined in $L^1\big( \varphi_i(I) \dd I \dd\bfv\big)$.
In particular, the number density $n_i$, mass density~$\varrho_i$,
mean velocity $\bfu_i$ and temperature $T_i$ of each species are respectively given by
\begin{equation}
n_i (t , {\bfx}) =  \int_{\RR^3} \! \int_0^{+\infty} \,
              f_i ( t, {\bfx} , {\bfv}, I ) \varphi_i(I) \,
             \dd I\, \dd{\bfv},
\label{eq:ni}
\end{equation}
\begin{equation}
\varrho_i (t , {\bfx}) = m_i n_i (t , {\bfx}) ,
\label{eq:ri}
\end{equation}
\begin{equation}
\bfu_i (t , {\bfx}) =  \frac{1}{n_i (t,\bfx)} \int_{\RR^3} \! \int_0^{+\infty} \,
               \bfv  f_i ( t, {\bfx} , {\bfv}, I ) \varphi_i(I) \,
             \dd I\, \dd{\bfv} ,
\label{eq:ui}
\end{equation}
\begin{equation}
T_i (t , {\bfx}) = \frac{1}{3 k_B n_i (t,\bfx)} \int_{\RR^3} \! \int_0^{+\infty} \,
              m_i  | \bfv - \bfu_i(t,\bfx) |^2  f_i ( t, {\bfx} , {\bfv}, I ) \varphi_i(I) \,
             \dd I\, \dd{\bfv}.
\label{eq:Ti}
\end{equation}


\subsection{Collisions and Borgnakke-Larsen procedure}
\label{ssec:cblp}

The particles of the mixture undergo binary collisions, either of elastic or reactive type.
Elastic collisions can occur among particles of the same constituent (mono-species elastic collisions) as well as among particles of different constituents (bi-species elastic collisions).
The mono-species and bi-species elastic collisions result in  changes in velocities and internal energies but do not modify the species and consequently do not modify the molecular masses of the colliding particles.
If we denote the velocities and internal energies of the colliding particles
before the collision by $\bfv_i$, $\bfv_j$ and $I_i$, $I_j$ respectively,
and their corresponding post-collisional values by $\bfv_i'$, $\bfv_j'$  and $I_i'$, $I_j'$,
the conservation laws of momentum and total energy for elastic collisions are given by
\begin{eqnarray}
& m_i \bfv_i + m_j \bfv_j = m_i \bfv_i' + m_j \bfv_j' ,
\label{eq:conservation_mom_elastic}
\\[1mm]
& \dps \frac12 m_i | \bfv_i |^2 + I_i + \frac12 m_j | \bfv_j |^2 + I_j =
      \frac12 m_i | \bfv_i' |^2 + I_i' + \frac12 m_j | \bfv_j' |^2 + I_j' .
\label{eq:conservation_eng_elastic}
\end{eqnarray}
In particular, $i\not=j$ for bi-species collisions
and $i=j$ for mono-species collisions.
In the latter case,  we use the indices $i$ and $i_*$ to distinguish the velocities {and internal energies} of the two  colliding particles.

\medskip

On the other hand,
reactive collisions occur among particles of constituents $A_1$, $A_2$ or $A_3$, $A_4$
and follow the reaction law \eqref{eq:cr}.
These reactive collisions result not just in changes in velocities and internal energies but also {result} in a transformation of the reactants into products of the reaction.
Consequently, they imply a rearrangement of mass and
a redistribution of chemical binding energy.
If  $A_i$, $A_j$ and $A_k$, $A_l$ represent the reactants and products of the chemical reaction, with
$(i,j,k,l) \in \{(1,2,3,4), (2,1,4,3)\}$  and $ (i,j,k,l) \in \left \{(3,4,1,2),(4,3,2,1) \right \} $ for the forward and  backward chemical reactions respectively,
the conservation laws of momentum and total energy (kinetic and internal energies as well as chemical binding energy)
for reactive collisions are given by
\begin{eqnarray}
& m_i \bfv_i + m_j \bfv_j = m_k \bfv_k' + m_l \bfv_l' ,
    \label{eq:conservation_mom_reactive}
    \\[1mm]
& \dps \frac12 m_i | \bfv_i |^2 + I_i + E_i + \frac12 m_j | \bfv_j |^2 + I_j + E_j  =
    \frac12 m_k | \bfv_k' |^2 + I_k' + E_k + \frac12 m_l | \bfv_l' |^2 + I_l' + E_l.
\label{eq:conservation_eng_chemical}
\end{eqnarray}
Equation \eqref{eq:conservation_eng_chemical} can be written in an equivalent form as
$$
\frac12 m_i | \bfv_i |^2 + I_i + \frac12 m_j | \bfv_j |^2 + I_j - \frac{E}{2}
   = \frac12 m_k | \bfv_k' |^2 + I_k' + \frac12 m_l | \bfv_l' |^2 + I_l' + \frac{E}{2},
    \eqno{(\ref{eq:conservation_eng_chemical}a)}
$$
where $E$ has been defined in~\eqref{eq:E}.


\bigskip

As usual in kinetic theory, the post-collisional velocities {can be} expressed in terms of pre-collisional velocities
and the corresponding expressions {for both elastic and reactive collisions} are derived from the conservation laws
\eqref{eq:conservation_mom_elastic},\eqref{eq:conservation_eng_elastic}
and  \eqref{eq:conservation_mom_reactive},\eqref{eq:conservation_eng_chemical} {respectively}.
In order to give these expressions, we describe first the Borgnakke-Larsen procedure \cite{bor-lar-75}.
Such procedure is based on the repartition of the total energy of the colliding pair
into kinetic and internal energies, when the collisions are of elastic type,
or into kinetic, internal and chemical binding energies, when the collisions are of reactive type.
Let us consider first elastic collisions and compute the total energy $\cal E$ of the colliding pair.
In the centre of mass reference frame, due to the conservation
laws \eqref{eq:conservation_mom_elastic} and \eqref{eq:conservation_eng_elastic},
we have that
\begin{eqnarray}
{\cal E} = \frac12 \mu_{ij} | \bfv_i - \bfv_j |^2 + I_i + I_j
             = \frac12 \mu_{ij} | \bfv_i' - \bfv_j' |^2 + I_i' + I_j' ,
\label{eq:cE}
\end{eqnarray}
where $\mu_{ij} = {m_i m_j}/{(m_i + m_j)}$ is the reduced mass of the colliding pair,
$\bfv_i - \bfv_j$ and $\bfv_i' - \bfv_j'$ are the relative velocities before and after the collision.
{Next} we introduce a parameter $R \!\in\! [0,1]$ and attribute the portions $R{\cal E}$ and $(1-R){\cal E}$ {of the total energy to the kinetic and internal energies,} respectively, of the outgoing pair,
that is
\begin{equation}
\frac12 \mu_{ij} | \bfv_i' - \bfv_j' |^2 = R{\cal E}
\qquad \mbox{and} \qquad
I_i' + I_j'  = (1-R) {\cal E} .
\label{eq:EKI}
\end{equation}
The first equation of \eqref{eq:EKI} can be parametrized by a unit vector $\omega\in \mathbb{S}^{2}$ to obtain
\begin{equation}
  \bfv_i' - \bfv_j'=\sqrt{\frac{2R{\cal E}}{\mu_{ij}}} \;
             T_\omega\left [ \frac{\bfv_i-{\bfv}_j}{\left |\bfv_i-{\bfv}_j \right |} \right ],
\label{eq:KE_parametrization}
\end{equation}
where
$T_{\omega} [ \bfx] = \bfx - 2({\omega}\cdot \bfx) {\omega}$
is the symmetry with respect to the plane
$\{\omega\}^\perp$.

Furthermore, we introduce another parameter, $r \!\in\! [0,1]$,
and allocate the portions $r(1-R){\cal E}$ and $(1-r)(1-R){\cal E}$
of the internal energy to each outgoing particle, that is
\begin{equation}
I_i' = r (1-R) {\cal E}
\qquad \mbox{and} \qquad
I_j' = (1-r) (1-R) {\cal E} .
\label{eq:IiIj}
\end{equation}

\medskip


\noindent
As a consequence of the Borgnakke-Larsen procedure \cite{DMS-EJM-05},
using the conservation law \eqref{eq:conservation_mom_elastic} together with~\eqref{eq:KE_parametrization},
we can express the elastic post-collisonal velocities
in terms of pre-collisonal velocities as given below
\begin{eqnarray}
& & {\bfv}'_i = \frac{m_i\bfv_i+m_j{\bfv}_j}{m_i+m_j} +
             \frac{m_j}{m_i+m_j} \sqrt{\frac{2R{\cal E}}{\mu_{ij}}} \;
             T_\omega\left [ \frac{\bfv_i-{\bfv}_j}{\left |\bfv_i-{\bfv}_j \right |} \right] ,
             \nonumber  \\
             \label{eq:elastic_bispecies_postcollisional_velocity}
             \\ [-2mm]
& & {\bfv}'_j=\frac{m_i\bfv_i+m_j{\bfv}_j}{m_i+m_j} -
             \frac{m_i}{m_i+m_j}\sqrt{\frac{2R{\cal E}}{\mu_{ij}}} \;
             T_\omega\left [ \frac{\bfv_i-{\bfv}_j}{\left |\bfv_i-{\bfv}_j \right |} \right] .
             \nonumber
\end{eqnarray}
In the particular case of mono-species elastic collisions, the post{-collisional} velocities are given by
\begin{eqnarray}
& & \bfv_i' = \frac{\bfv_i+{\bfv}_{i_*}}2 +
             \sqrt{\frac{R{\cal E}}{m_i}} \;
             T_\omega\left [ \frac{\bfv_i-{\bfv}_{i_*}}{\left |\bfv_i-{\bfv}_{i_*} \right |} \right] ,
             \nonumber  \\
             \\[-2mm]
& & {\bfv}'_{i_*}=\frac{{\bfv}_{i}+{\bfv}_{i_*}}{2} -
             \sqrt{\frac{R{\cal E}}{m_i }} \;
             T_\omega\left [ \frac{\bfv_i-{\bfv}_{i_*}}{\left |\bfv_i-{\bfv}_{i_*} \right |} \right] .
             \nonumber
\end{eqnarray}

\medskip


\noindent
Concerning now reactive collisions associated to the forward chemical reaction,
we denote by ${\cal E}^*$ the total energy of the colliding pair.
In the centre of mass reference frame, due to the conservation laws \eqref{eq:conservation_mom_reactive} and~(\ref{eq:conservation_eng_chemical}$a$),
we have that
\begin{eqnarray}
{\cal E}^* = \frac12 \mu_{12} | \bfv_1 - \bfv_2 |^2 + I_1 + I_2 - \frac{E}{2}
                = \frac12 \mu_{34} | \bfv_3' - \bfv_4' |^2 + I_3' + I_4' + \frac{E}{2} .
\label{eq:Ef}
\end{eqnarray}
Introducing the parameters $R, r \!\in\! [0,1]$ and assuming the following repartition of the {total energy into} kinetic  and internal energies of the  outgoing particles, namely
\begin{equation}
\frac12 \mu_{34} | \bfv_3' - \bfv_4' |^2 = R{\cal E}^* - \frac{E}{6} ,
\qquad
I_3' + I_4'  = (1-R) {\cal E}^* - \frac{E}{3} ,
\label{eq:EKIf}
\end{equation}
with
\begin{equation}
I_3' = r (1-R) {\cal E}^* - \frac{E}{6} ,
\qquad
I_4' = (1-r) (1-R) {\cal E}^* - \frac{E}{6} ,
\label{eq:Ei34}
\end{equation}
we obtain that the first equation of \eqref{eq:EKIf} can be parameterized by a unit vector $\omega\in \mathbb{S}^{2}$
as given below
\begin{equation}
\bfv_3' - \bfv_4'= {\sqrt{\frac{2}{\mu_{34}}} \left (R{\cal E}^*-\frac{E}{6} \right )^{1/2}} \;
             T_\omega\left [ \frac{{\bfv}_1-{\bfv}_2}{\left |{\bfv}_1-{\bfv}_2 \right |} \right ].
\label{eq:KE_parametrization1}
\end{equation}
Using \eqref{eq:KE_parametrization1} together with the conservation law \eqref{eq:conservation_mom_reactive},
we can express the reactive post-collisional velocities for the forward reaction in terms of  pre-collisional velocities as given below
\begin{equation}
\begin{aligned}
{\bfv}_3'  &= \frac{m_1{\bfv}_1+m_2{\bfv}_2}{m_1+m_2}{+}\frac{m_4}{m_3+m_4}
                    \sqrt{\frac{2}{\mu_{34}}} \left (R{\cal E}^*-\frac{E}{6} \right )^{1/2}
                    T_\omega \left [ \frac{{\bfv}_1-{\bfv}_2}{\left | {\bfv}_1-{\bfv}_2\right |} \right] ,
                    \\[2mm]
{\bfv}_4'   &= \frac{m_1{\bfv}_1+m_2{\bfv}_2}{m_1+m_2}{-}\frac{m_3}{m_3+m_4}
                    \sqrt{\frac{2}{\mu_{34}}}  \left (R{\cal E}^*-\frac{E}{6} \right )^{1/2}
                    T_\omega \left [ \frac{{\bfv}_1-{\bfv}_2}{\left | {\bfv}_1-{\bfv}_2\right |} \right] .
\end{aligned}
\label{eq:postcollisional_velocity_forward_reaction}
\end{equation}
Analogously, for reactive collisions associated to the backward chemical reaction,
we write the conservation law of total energy in the centre of mass reference frame
in the form
\begin{eqnarray}
{\cal E}^* = \frac12 \mu_{34} | \bfv_3 - \bfv_4 |^2 + I_3 + I_4 + \frac{E}{2}
               = \frac12 \mu_{12} | \bfv_1' - \bfv_2' |^2 + I_1' + I_2' - \frac{E}{2} .
\label{eq:Eb}
\end{eqnarray}
Then, as before, we assume {the following repartition of the  total energy into kinetic
and internal energies of the  outgoing particles, namely}
\begin{equation}
\frac12 \mu_{12} | \bfv_1' - \bfv_2' |^2 = R{\cal E}^* + \frac{E}{6} ,
\qquad
I_1' + I_2'  = (1-R) {\cal E}^* + \frac{E}{3} ,
\label{eq:EKIb}
\end{equation}
with
\begin{equation}
I_1' = r (1-R) {\cal E}^* + \frac{E}{6} ,
\qquad
I_2' = (1-r) (1-R) {\cal E}^* + \frac{E}{6} .
\label{eq:Ei12}
\end{equation}
The first equation of \eqref{eq:EKIb} can be parameterized by a unit vector $\omega\in \mathbb{S}^{2}$ as given below
\begin{equation}
\bfv_1' - \bfv_2'= \sqrt{\frac{2}{\mu_{12}}} \left (R{\cal E}^*+\frac{E}{6} \right )^{1/2} \;
                           T_\omega\left [ \frac{{\bfv}_3-{\bfv}_4}{\left |{\bfv}_3-{\bfv}_4 \right |} \right ].
\label{eq:KE_parametrization2}
\end{equation}
Using \eqref{eq:KE_parametrization2} together with the conservation law \eqref{eq:conservation_mom_reactive},
we can express the reactive post-collisional velocities for the backward reaction in terms of  pre-collisional velocities as given below
\begin{equation}
\begin{aligned}
{\bfv}_1'  &= \frac{m_3{\bfv}_3+m_4{\bfv}_4}{m_3+m_4}{+}\frac{m_2}{m_1+m_2}
                    \sqrt{\frac{2}{\mu_{12}}} \left (R{\cal E}^* + \frac{E}{6} \right )^{1/2}
                    T_\omega \left [ \frac{{\bfv}_3-{\bfv}_4}{\left | {\bfv}_3-{\bfv}_4\right |} \right] ,
                    \\[2mm]
{\bfv}_2'   &= \frac{m_3{\bfv}_3+m_4{\bfv}_4}{m_3+m_4}{-}\frac{m_1}{m_1+m_2}
                    \sqrt{\frac{2}{\mu_{12}}}  \left (R{\cal E}^* + \frac{E}{6} \right )^{1/2}
                    T_\omega \left [ \frac{{\bfv}_3-{\bfv}_4}{\left | {\bfv}_3-{\bfv}_4\right |} \right] .
\end{aligned}
\label{eq:postcollisiona_velocity_reverse_reaction}
\end{equation}
For other details on the collisional dynamics, the reader is referred to \cite{DMS-EJM-05}.


\subsection{Kinetic equations}
\label{ssec:keqs}

The time-space evolution of the distribution functions {$f_i(\bfv,I)$, with $i=1,\ldots,4$,}
is specified by the system of kinetic equations
\begin{equation}
\frac{\partial f_i}{\partial t} + \bfv \cdot\nabla_x f_i = \sum_{j=1}^4Q_{ij}^e(f_i,f_j) + Q_{i}^{react},
\qquad i=1,\ldots,4 .
\label{eq:ke}
\end{equation}
Above, {when $i\not=j$, the notation $Q_{ij}^e$ represents the bi-species elastic} operator associated to collisions between one particle of constituent $i$ and another one of constituent $j$, whereas when $i=j$, it reduces to the mono-species elastic collisional operator $Q_{ii}^e$.
Moreover, $Q_{i}^{react}$ represents the reactive collisional operator.
The operators $Q_{ij}^e$ and $Q_{i}^{react}$ are defined as follows.
For bi-species non-reactive interactions (i.e. $i\not=j)$, the operators take the form
\begin{eqnarray}
Q_{ij}^e(f_i,f_j) &=&  \int_{\RR^3} \! \int_0^{+\infty} \! \int_0^1 \! \int_0^1 \! \int_{\mathbb{S}^{2}} \!\!
             \left [f_i(\bfv_i', I_i') f_j({\bfv}_j', I'_j) \!-\! f_i(\bfv, I) f_j({\bfv}_j, I_j)   \right ]
             \label{eq:bispecies} \\
             & & \hspace*{0.5cm}
             \times B_{ij} (\bfv,{\bfv}_j,I,I_j,R,r,\omega) (1-R)\left | \bfv-{\bfv}_j \right |^{-1} \varphi_i(I)^{-1}\,
             \dd\omega\,\dd r\, \dd R\, \dd I_j\, \dd{\bfv}_j ,
             \hspace*{0.5cm}
             \nonumber
\end{eqnarray}
where  $B_{ij}$ are suitable cross sections and ${\bfv}'_i$, ${\bfv}'_j$, $I_i'$ and $I_j'$ are given by
\begin{eqnarray*}
& & {\bfv}'_i = \frac{m_i\bfv+m_j{\bfv}_j}{m_i+m_j} +
             \frac{m_j}{m_i+m_j} \sqrt{\frac{2R{\cal E}}{\mu_{ij}}} \;
             T_\omega\left [ \frac{\bfv-{\bfv}_j}{\left |\bfv-{\bfv}_j \right |} \right] ,
             \nonumber  \\
             \\ [-2mm]
& & {\bfv}'_j=\frac{m_i\bfv+m_j{\bfv}_j}{m_i+m_j} -
             \frac{m_i}{m_i+m_j}\sqrt{\frac{2R{\cal E}}{\mu_{ij}}} \;
             T_\omega\left [ \frac{\bfv-{\bfv}_j}{\left |\bfv-{\bfv}_j \right |} \right],
             \nonumber
\end{eqnarray*}
\begin{equation*}
I_i' = \left({\frac12} \mu_{ij} | \bfv - \bfv_j |^2 + I + I_j\right )   r (1-R)
\qquad \mbox{and} \qquad
I_j' = \left( {\frac12} \mu_{ij} | \bfv - \bfv_j |^2 + I + I_j\right ) (1-r) (1-R).
\end{equation*}
In the case of mono-species elastic collisions, the operators have the following structure
\begin{eqnarray}
Q_{ii}^e(f_i,f_{i}) &=&  \int_{\RR^3} \! \int_0^{+\infty} \! \int_0^1 \! \int_0^1 \! \int_{\mathbb{S}^{2}} \!\!
             \left [f_i(\bfv_i', I_i') f_i({\bfv}'_{i_*}, I'_{i_*}) \!-\! f_i(\bfv, I) f_i({\bfv}_{i_*}, I_{i_*})   \right ]
             \\
             & & \hspace*{0.5cm}
             \times B_{ii} (\bfv,{\bfv}_{i_*},I,I_{i_*},R,r,\omega) (1-R)\left | \bfv-{\bfv}_{i_*} \right |^{-1} \varphi_i(I)^{-1}\,
             \dd\omega\,\dd r\, \dd R\, \dd I_{i_*}\, \dd{\bfv}_{i_*} ,
             \hspace*{0.5cm}
             \nonumber
\end{eqnarray}
where  $B_{ii}$ are suitable cross sections and ${\bfv}'_i$, ${\bfv}'_{i_*}$, $I_i'$ and $I_{i_*}$ are given by
\begin{eqnarray*}
& & \bfv_i' = \frac{\bfv+{\bfv}_{i_*}}2 +
             \sqrt{\frac{R{\cal E}}{m_i}} \;
             T_\omega\left [ \frac{\bfv-{\bfv}_{i_*}}{\left |\bfv-{\bfv}_{i_*} \right |} \right] ,
             \nonumber  \\
             \label{eq:elastic_monospecies_postcollisional_velocity}
             \\[-2mm]
& & {\bfv}'_{i_*}=\frac{{\bfv}+{\bfv}_{i_*}}{2} -
             \sqrt{\frac{R{\cal E}}{m_i }} \;
             T_\omega\left [ \frac{\bfv-{\bfv}_{i_*}}{\left |\bfv-{\bfv}_{i_*} \right |} \right],
             \nonumber
\end{eqnarray*}
\begin{equation*}
I_i' = \left({\frac{m_i}{4} }| \bfv - \bfv_{i_*} |^2 + I + I_{i_*}\right )   r (1-R)
\quad \mbox{and} \quad
I_{i_*}' = \left({\frac{m_i}{4} } | \bfv - \bfv_{i_*} |^2 + I + I_{i_*}\right ) (1-r) (1-R).
\end{equation*}
The reactive collisional operators are defined in a more involved way. 
Consider for some suitable sets $F_i$ the Heaviside-like function
\begin{equation*}
\dps H_i(\xi) \!=\! \left\{
                        \begin{array}{ll}
                        1 ,  \quad & \xi \in F_i \\[0.5mm]
                        0 , & \xi\notin F_i,
                        \end{array}
                        \right.
                        \qquad i=1,\dots, 4,
\end{equation*}
and let 
$B^{react}\, :\, \mathbb{R}^{3}\times\mathbb{R}^{3} \times \mathbb{R}^+\times \mathbb{R}^+\times [0,1]^2\times  \mathbb{S}^{2} \to \mathbb{R}^+$ 
be a suitable cross-section.

\noindent
We first treat the case of the forward reaction.
{To do this, we define the total energy}
\begin{eqnarray*}
{\cal E}_1^* = \frac12 \mu_{12} | \bfv - \bfv_2 |^2 + I + I_2 - \frac{E}{2},
\end{eqnarray*}
{the reactive post-collisional velocities} and internal energies
\begin{equation*}
\begin{aligned}
{\bfv}_3'  &= \frac{m_1{\bfv}+m_2{\bfv}_2}{m_1+m_2}{+}\frac{m_4}{m_3+m_4}
                    \sqrt{\frac{2}{\mu_{34}}} \left (R{\cal E}_1^*-\frac{E}{6} \right )^{1/2}
                    T_\omega \left [ \frac{{\bfv}-{\bfv}_2}{\left | {\bfv}-{\bfv}_2\right |} \right] ,
                    \\[2mm]
{\bfv}_4'   &= \frac{m_1{\bfv}+m_2{\bfv}_2}{m_1+m_2}{-}\frac{m_3}{m_3+m_4}
                    \sqrt{\frac{2}{\mu_{34}}}  \left (R{\cal E}_1^*-\frac{E}{6} \right )^{1/2}
                    T_\omega \left [ \frac{{\bfv}-{\bfv}_2}{\left | {\bfv}-{\bfv}_2\right |} \right] ,
                     \\[2mm]
I_3' &= r (1-R) {\cal E}_1^* - \frac{E}{6} ,
 	             \\[2mm]
I_4' &= (1-r) (1-R) {\cal E}_1^* - \frac{E}{6} ,
\end{aligned}
\end{equation*}
{as well as} the admissible set
\begin{eqnarray*}
\lefteqn{
F_1 \!=\! \bigg\{
             (I , I_2 , r , R, \bfv , {\bfv}_2 ) : \;\;
             I\geq \frac{E}{6},  \; I_2\geq \frac{E}{6}, \;
             R{\cal E}_1^*\geq \frac{E}{6}, }
             \\[-1mm]
             & & \hspace*{3.75cm}  \frac{\mu_{12}}{2}\left | \bfv \!-\! {\bfv}_2 \right |^2\geq \frac{E}{6}, \;
             (1\!-\!R)r{\cal E}_1^*\geq \frac{E}{6}, \;
             (1\!-\!R)(1\!-\!r){\cal E}_1^*\geq \frac{E}{6}
             \bigg\} .
\end{eqnarray*}
The first collisional integral describing the forward chemical reaction is hence
\begin{multline}
Q_{1}^{react} (\bfv, I)= \int_{\mathbb{R}^3}\int_0^{+\infty} \int_0^1\int_0^1\int_{\mathbb{S}^{2}}  \!
            \left [ \! \left (\frac{m_1 m_2}{m_3 m_4} \right )^{3} \!\!
            f_3({\bfv}'_3, I'_3)f_4({\bfv}'_4, I'_4) - f_1(\bfv, I) f_2({\bfv}_2, I_2) \right ]
            \\
            \times H_1(I , I_2 , r , R, \bfv , {\bfv}_2 )B^{react}(\bfv,{\bfv}_2, I,I_2,R,r,\omega) \,
            (m_1m_2)^{-2}\left | \bfv-{\bfv}_2\right |^{-1}
            \\
            \times
            (1-R) \varphi_1(I)^{-1}\, \dd\omega\, \dd r\, \dd R\, \dd I_2\, \dd{\bfv}_2 .
\label{Q:reactive_1}
\end{multline}
The structure of $Q_{2}^{react} $ is similar. We define {the total energy}
\begin{eqnarray*}
{\cal E}_2^* = \frac12 \mu_{12} | \bfv - \bfv_1 |^2 + I + I_1 - \frac{E}{2},
\end{eqnarray*}
the reactive post-collisional velocities and internal energies
\begin{equation*}
\begin{aligned}
{\bfv}_3'  &= \frac{m_1{\bfv_1}+m_2{\bfv}}{m_1+m_2}{+}\frac{m_4}{m_3+m_4}
                    \sqrt{\frac{2}{\mu_{34}}} \left (R{\cal E}_2^*-\frac{E}{6} \right )^{1/2}
                    T_\omega \left [ \frac{{\bfv_1}-{\bfv}}{\left | {\bfv_1}-{\bfv}\right |} \right] ,
                    \\[2mm]
{\bfv}_4'   &= \frac{m_1{\bfv}_1+m_2{\bfv}}{m_1+m_2}{-}\frac{m_3}{m_3+m_4}
                    \sqrt{\frac{2}{\mu_{34}}}  \left (R{\cal E}_2^*-\frac{E}{6} \right )^{1/2}
                    T_\omega \left [ \frac{{\bfv}_1-{\bfv}}{\left | {\bfv}_1-{\bfv}\right |} \right] ,
                     \\[2mm]
I_3' &= r (1-R) {\cal E}_2^* - \frac{E}{6} ,
 	             \\[2mm]
I_4' &= (1-r) (1-R) {\cal E}_2^* - \frac{E}{6} ,
\end{aligned}
\end{equation*}
{as well as} the admissible set
\begin{eqnarray*}
\lefteqn{
F_2 \!=\! \bigg\{
             (I , I_1 , r , R, \bfv , {\bfv}_1 ) : \;\;
             I\geq \frac{E}{6},  \; I_1\geq \frac{E}{6}, \;
             R{\cal E}_2^*\geq \frac{E}{6}, }
             \\
             & & \hspace*{3.75cm}
             \frac{\mu_{12}}{2}\left | {\bfv \!-\! {\bfv}_1} \right |^2\geq \frac{E}{6}, \;
             (1\!-\!R)r{\cal E}_2^*\geq \frac{E}{6}, \;
             (1\!-\!R)(1\!-\!r){\cal E}_2^*\geq \frac{E}{6}
             \bigg\} .
\end{eqnarray*}
The collisional integral $Q_{2}^{react}$ describing the forward chemical reaction is hence
\begin{multline}
Q_{2}^{react} (\bfv, I)= \int_{\mathbb{R}^3}\int_0^{+\infty} \int_0^1\int_0^1\int_{\mathbb{S}^{2}}  \!
            \left [ \! \left (\frac{m_1 m_2}{m_3 m_4} \right )^{3} \!\!
            f_3({\bfv}'_3, I'_3) f_4({\bfv}'_4, I'_4) - f_2(\bfv, I) f_1({\bfv}_1, I_1) \right ]
            \\
            \times H_2(I , I_1 , r , R, \bfv , {\bfv}_1 )B^{react}(\bfv,{\bfv}_1, I,I_1,R,r,\omega) \,
            (m_1m_2)^{-2}\left | {\bfv-{\bfv}_1}\right |^{-1}
            \\
            \times
            (1-R) \varphi_2(I)^{-1}\, \dd\omega\, \dd r\, \dd R\, \dd I_1\, \dd{\bfv}_1 .
\label{Q:reactive_2}
\end{multline}
In the case of the backward reaction, we have to treat two cases. Let
\begin{eqnarray*}
{\cal E}_3^* = \frac12 \mu_{34} | \bfv - \bfv_4 |^2 + I + I_4 + \frac{E}{2}
\end{eqnarray*}
{be the total energy.} 
Then, as before, we write the {reactive post collisional velocities} and internal energies, respectively as given below
\begin{equation*}
\begin{aligned}
{\bfv}_1'  &= \frac{m_3{\bfv}+m_4{\bfv}_4}{m_3+m_4}{+}\frac{m_2}{m_1+m_2}
                    \sqrt{\frac{2}{\mu_{12}}} \left (R{\cal E}_3^* + \frac{E}{6} \right )^{1/2}
                    T_\omega \left [ \frac{{\bfv}-{\bfv}_4}{\left | {\bfv}-{\bfv}_4\right |} \right] ,
                    \\[2mm]
{\bfv}_2'   &= \frac{m_3{\bfv}+m_4{\bfv}_4}{m_3+m_4}{-}\frac{m_1}{m_1+m_2}
                    \sqrt{\frac{2}{\mu_{12}}}  \left (R{\cal E}_3^* + \frac{E}{6} \right )^{1/2}
                    T_\omega \left [ \frac{{\bfv}-{\bfv}_4}{\left | {\bfv}-{\bfv}_4\right |} \right] ,
                    \\[2mm]
I_1' &= r (1-R) {\cal E}_3^* + \frac{E}{6} ,
                    \\[2mm]
I_2' &= (1-r) (1-R) {\cal E}_3^* + \frac{E}{6} .
\end{aligned}
\end{equation*}
The admissible set $F_3$ is
$$
F_3 \!=\!  \Big\{ (I_4 , r , R , {\bfv}_4 ) : \;\;   I_4\geq 0, \;\; {\bfv}_4\in\mathbb{R}^3, \;\; r, R \in [0,1] \Big\},
$$
and hence the reactive {collisional integral $Q_{3}^{react}$  describing the backward reaction} is defined by
\begin{multline}
Q_{3}^{react} (\bfv, I)= \int_{\mathbb{R}^3}\int_0^{+\infty} \int_0^1\int_0^1\int_{\mathbb{S}^{2}}  \!
            \left [ \! \left (\frac{m_3 m_4}{m_1 m_2} \right )^{3} \!\!
            f_1({\bfv}'_1, I'_1)f_2({\bfv}'_2, I'_2) - f_3(\bfv, I) f_4({\bfv}_4, I_4 )\right ]
            \\[1mm]
            \times
            H_3(I_4 , r , R , {\bfv}_4) B^{react}(\bfv,{\bfv}_4, I,I_4,R,r,\omega) \,
            (m_3m_4)^{-2}\left | \bfv-{\bfv}_4\right |^{-1}
            \\[1mm]
            \times (1-R) \varphi_3(I)^{-1}\, \dd\omega\, \dd r\, \dd R\, \dd I_4\, \dd{\bfv}_4 .
\label{Q:reactive_3}
\end{multline}
We conclude the description of the reactive collisional operators by {defining $Q_{4}^{react}$. Let}
\begin{eqnarray*}
{\cal E}_4^* = \frac12 \mu_{34} | \bfv - \bfv_3 |^2 + I + I_3 + \frac{E}{2}
\end{eqnarray*}
{be the total energy.} 
Then, as before, we write the {reactive post collisional} velocities and internal energies, respectively as given below
\begin{equation*}
\begin{aligned}
{\bfv}_1'  &= \frac{m_3{\bfv}_3+m_4{\bfv}}{m_3+m_4}{+}\frac{m_2}{m_1+m_2}
                    \sqrt{\frac{2}{\mu_{12}}} \left (R{{\cal E}_4^*} + \frac{E}{6} \right )^{1/2}
                    T_\omega \left [ \frac{{\bfv}_3-{\bfv}}{\left | {\bfv}_3-{\bfv}\right |} \right] ,
                    \\[2mm]
{\bfv}_2'   &= \frac{m_3{\bfv}_3+m_4{\bfv}}{m_3+m_4}{-}\frac{m_1}{m_1+m_2}
                    \sqrt{\frac{2}{\mu_{12}}}  \left (R{{\cal E}_4^*} + \frac{E}{6} \right )^{1/2}
                    T_\omega \left [ \frac{{\bfv}_3-{\bfv}}{\left | {{\bfv}_3-{\bfv}}\right |} \right] .
                    \\[2mm]
I_1' &= r (1-R) {\cal E}_4^* + \frac{E}{6} ,
                    \\[2mm]
I_2' &= (1-r) (1-R) {\cal E}_4^* + \frac{E}{6} .
\end{aligned}
\end{equation*}
The admissible set $F_4$ is
$$
F_4 \!=\!  \Big\{ (I_3 , r , R , {\bfv}_3 ) : \;\;   I_3\geq 0, \;\; {\bfv}_3\in\mathbb{R}^3, \;\; r, R \in [0,1] \Big\},
$$
and hence the reactive {collisional integral $Q_{4}^{react}$ describing the backward reaction} is defined by
\begin{multline}
Q_{4}^{react} (\bfv, I)= \int_{\mathbb{R}^3}\int_0^{+\infty} \int_0^1\int_0^1\int_{\mathbb{S}^{2}}  \!
            \left [ \! \left (\frac{m_3 m_4}{m_1 m_2} \right )^{3} \!\!
            f_1({\bfv}'_1, I'_1)f_2({\bfv}'_2, I'_2) - f_4(\bfv, I) f_3({\bfv}_3, I_3 )\right ]
            \\[2mm]
            \times H_4(I_3 , r , R , {\bfv}_3)B^{react}(\bfv,{\bfv}_3, I,I_3,R,r,\omega) \,
            (m_3m_4)^{-2}\left | {\bfv-{\bfv}_3}\right |^{-1}
            \\[2mm]
           \times (1-R) \varphi_4(I)^{-1}\, \dd\omega\, \dd r\, \dd R\, \dd I_3\, \dd{\bfv}_3 .
\label{Q:reactive_4}
\end{multline}


\subsection{Properties of the collisional operators}
\label{ssec:pco}
Here we review some properties of the collisional operators given in \cite{DMS-EJM-05}
that will be used in the derivation of the limit equations.
\begin{lem}
\label{lem:Q_im}
(See Lemma $2$, page $223$ of {\it \cite{DMS-EJM-05}})
\\
Let ${\psi} \! : \mathbb{R}^3\times[0,+\infty)\rightarrow \mathbb{R}$ be a function such that the weak formulation
\begin{equation*}
\int_{\mathbb{R}^3}\int_0^{+\infty} {Q_{ii}^e} \psi(\bfv,I) \varphi_i(I) \dd I \dd\bfv
\end{equation*}
makes sense.     
Then 
\begin{eqnarray}
\lefteqn{\int_{\mathbb{R}^3}\int_0^{+\infty} {Q_{ii}^e} \psi(\bfv,I) \varphi_i(I) \dd I \dd\bfv }
               \nonumber  \\
      &=&  - \frac{1}{4}\int_{\mathbb{R}^3}\int_{\mathbb{R}^3}\int_0^{+\infty}\int_0^{+\infty}
               \int_0^1\int_0^1\int_{\mathbb{S}^{2}}
               {\left [f_i(\bfv_i', I_i') f_i({\bfv}'_{i_*}, I'_{i_*}) \!-\! f_i(\bfv, I) f_i({\bfv}_{i_*}, I_{i_*})   \right ]}
               \label{eq:weakform_mono}   \\[1mm]
      &  &  \times \, \Big[ \psi(\bfv_i',I_i')+\psi({\bfv}'_{i_*},I'_{i_*})-\psi(\bfv,I)-\psi({\bfv}_{i_*},I_{i_*}) \Big ]
               {B_{ii}}(\bfv,{\bfv}_{i_*},I,I_{i_*},R,r,\omega)
               \nonumber   \\[1mm]
      &  &  \times \, (1-R)\left | \bfv-{\bfv}_{i_*} \right |^{-1}\,\dd\omega\,\dd r\, \dd R\, \dd I_{i_*}\,\dd I\, \dd{\bfv}_{i_*}\,\dd\bfv.
               \nonumber
\end{eqnarray}
\end{lem}
\begin{lem}
\label{lem:Q_iplusQ_j}
(See Lemma $3$, page $224$ of {\it\cite{DMS-EJM-05}})
\\
Let $j\not =i$ and ${\psi} \! : \mathbb{R}^3\times[0,+\infty)\rightarrow \mathbb{R}$ be a function such that the formulas
\begin{equation*}
\int_{\mathbb{R}^3}\int_0^{+\infty} Q_{ij}^e \psi(\bfv,I) \varphi_i(I) \dd I \dd\bfv
\qquad \text{and} \qquad
\int_{\mathbb{R}^3}\int_0^{+\infty} Q_{ji}^e \psi{(\bfv_j,I_j) \varphi_j(I_j) \dd I_j \dd\bfv_j}
\end{equation*}
make sense. Then
\begin{eqnarray}
\lefteqn{\int_{\mathbb{R}^3}\int_0^{+\infty} Q_{ij}^e\psi(\bfv,I)\varphi_i(I) \dd I \dd\bfv}
             \label{eq:weakform_bispecies_collision_integral1}   \\
            &\!\!=\!\!& -\frac{1}{2}\int_{\mathbb{R}^3}\!\!\int_{\mathbb{R}^3}\!\!\int_0^{+\infty}\!\!\!
                    \int_0^{+\infty}\!\!\!\int_0^1\!\int_0^1\!\int_{\mathbb{S}^{2}}\!
                   {\left [f_i(\bfv_i', I_i') f_j({\bfv}_j', I'_j) \!-\! f_i(\bfv, I) f_j({\bfv}_j, I_j)   \right ]}
                    \Big[ \psi(\bfv_i',I_i')-\psi(\bfv,I) \Big]
                    \nonumber \\
            &  & \times B_{ij}(\bfv,{\bfv}_j,I,I_j,R,r,\omega)(1-R)\left | \bfv-{\bfv}_j \right |^{-1}\,
                    \dd\omega\,\dd r\, \dd R\, \dd I_j\,\dd I \dd{\bfv}_j \, \dd\bfv
                    \nonumber
\end{eqnarray}
and
\begin{eqnarray}
\lefteqn{\int_{\mathbb{R}^3}\int_0^{+\infty} Q_{ij}^e \psi(\bfv,I) \varphi_i(I) \dd I \dd\bfv
            + \int_{\mathbb{R}^3}\int_0^{+\infty} Q_{ji}^e \psi({\bfv}_j,I_j)\varphi_j(I_j) \dd I_j \dd{\bfv}_j}
            \label{eq:weakform_bispecies_collision_integral2} \\
            &\!\!=\!\!& -\frac{1}{2}\int_{\mathbb{R}^3}\int_{\mathbb{R}^3}\int_0^{+\infty}\!\!\!
                    \int_0^{+\infty}\!\!\int_0^1\int_0^1{\int_{\mathbb{S}^{2}}} \!
                    {\left [f_i(\bfv_i', I_i') f_j({\bfv}_j', I'_j) \!-\! f_i(\bfv, I) f_j({\bfv}_j, I_j)   \right ]}
                    \nonumber  \\
            &  &  \times \Big[ \psi(\bfv_i',I_i')+\psi({\bfv}_j',I_j')-
                    \psi(\bfv,I) - \psi({\bfv}_j,I_j) \Big]B_{ij}(\bfv,{\bfv}_j,I,I_j,R,r,\omega)
                    \nonumber  \\
            &  &  \times
                    (1-R)\left | \bfv-{\bfv}_j \right |^{-1}\,\dd\omega\,\dd r\, \dd R\, \dd I_j\,\dd I \dd{\bfv}_j \, \dd\bfv.
                    \nonumber
\end{eqnarray}
\end{lem}
\begin{lem}
\label{lm:react}
(See Lemma $5$, page $229$ of {\it\cite{DMS-EJM-05}})
\\
Let ${\psi} \! :\mathbb{R}^3\times[0,+\infty)\rightarrow \mathbb{R}$
be a function such that for all $i=1,2,3,4,$ the formula
\begin{equation*}
\int_{\mathbb{R}^3}\int_0^{+\infty} Q_{i}^{react} \psi(\bfv,I) \varphi_i(I) \dd I \dd\bfv
\end{equation*}
makes sense.
Additionally, if $B^{react}({\color{black}{\bfv}_i},{\bfv}_j, {\color{black}I_i},I_j,R,r,\omega)$ for  $(i,j,k,l) \in \left\{ (1,2,3,4),(2,1,4,3) \right\}$
is equal to $B^{react}({\color{black}{\bfv}_i},{\bfv}_j, {\color{black}I_i},I_j,R,r,\omega)$ for $(i,j,k,l) \in \left\{ (3,4,1,2),(4,3,2,1) \right\} $,
then
\begin{eqnarray}
\lefteqn{\sum_{i=1}^{4}\int_{\mathbb{R}^3}\int_0^{+\infty} Q_{i}^{react} \psi(\bfv,I) \varphi_i(I) \dd I \dd\bfv}
            \label{eq:weakform_reactive_collision_integral} \\[1mm]
            &\!\!=\!\!& - \int_{\mathbb{R}^3} \int_{\mathbb{R}^3}
                            \int_0^{+\infty} \int_0^{+\infty} \int_0^1 \int_0^1 {\int_{\mathbb{S}^2}}
                            {\left[ {\frac{f_k(\bfv_k', I_k') f_l({\bfv}_l', I'_l)}{(m_k m_l)^3}} -
                            {\frac{f_i(\bfv_i, I_i) f_j({\bfv}_j, I_j)}{(m_i m_j)^3}} \right]}  \!
                            H_i
                            \nonumber \\[1mm]
            &  & \times \Big[ \psi({\bfv}'_k, I'_k)+\psi({\bfv}'_l, I'_l)-\psi({\bfv_i, I_i})-\psi({\bfv}_j, I_j) \Big]                
                            \nonumber  \\[1mm]
            &  &  \times B^{react}({\color{black}{\bfv}_i},{\bfv}_j, {\color{black}I_i},I_j,R,r,\omega)
                     m_im_j(1-R)\left | {\bfv_i}-{\bfv}_j \right |^{-1}\dd\omega\,\dd r\, \dd R\, \dd I_j\,{\dd I_i} \dd{\bfv}_j \, {\dd\bfv_i}.
                            \nonumber
\end{eqnarray}
\end{lem}


\subsection{Conservation laws and chemical rates}
\label{ssec:cl}

The conservation laws of the model are obtained from the properties stated in Lemmas
\ref{lem:Q_im}, \ref{lem:Q_iplusQ_j} and~\ref{lm:react} of Subsection \ref{ssec:pco}.


\begin{lem}
\label{lm:cls}
(See equations $(22)$, $(27)$, $(28)$, $(45)$, pages $224$, $225$, $226$, $231$ of $\cite{DMS-EJM-05}$)
\\
Consider the functions $ \Psi (\bfv , I)$ defined by
$ \Psi \!=\! (1,0,1,0)$, $\Psi \!=\! (1,0,0,1)$, $\Psi \!=\! (0,1,1,0)$
and also $ \Psi (\bfv , I)$ such that
$ \Psi \!=\! (\psi_1, \psi_2, \psi_3, \psi_4)$,
with $\psi_i = m_i \bfv_{x_1}$ or $m_i \bfv_{x_2}$ or $m_i \bfv_{x_3}$,
as well as
$ \Psi \!=\! (\psi_1, \psi_2, \psi_3, \psi_4)$, with $\psi_i = \frac12 m_i \bfv^2 + I + E_i $.
Then,
\begin{equation}
\sum_{i=1}^4 \! {\int_{\RR^3}} \! \int_0^{+\infty}
            \bigg( \sum_{j=1}^4 Q_{ij}^e + Q_i^{react} \bigg)
            \psi_i(\bfv , I) \,
            \varphi_i(I) \, \dd I\, \dd\bfv =  0.
\label{eq:cls}
\end{equation}
\QEDB
\end{lem}


\noindent
Property \eqref{eq:cls} indicates that elastic and reactive collisional operators are consistent with conservation
of physical quantities, namely
partial number densities $n_1\!+\!n_3$, $ n_1\!+\!n_4$, $n_2\!+\!n_3$,
momentum and total energy of the whole mixture
(kinetic, internal and chemical binding).


\begin{lem}
\label{lm:cpr}
The reactive collisional operators are such that
\begin{eqnarray}
\int_{\RR^3} \! \int_0^{+\infty}  \!\! Q_1^{react} \,  \varphi_1(I_1) \, \dd I_1 \, \dd\bfv_1
       \!&\!\!=\!\!&\!
       \int_{\RR^3} \! \int_0^{+\infty}  \!\! Q_2^{react} \,  \varphi_2(I_2) \, \dd I_2 \, \dd\bfv_2
       \label{eq:cpr}
       \\
       \!\!&\!\!=\!\!&\!\!
       - \int_{\RR^3} \!\int_0^{+\infty}  \!\! Q_3^{react} \, \varphi_3(I_3) \, \dd I_3 \, \dd\bfv_3
       \nonumber
       \\
       \!\!&\!\!=\!\!&\!\!
       - \int_{\RR^3} \! \int_0^{+\infty}  \!\! Q_4^{react} \, \varphi_4(I_4) \, \dd I_4 \, \dd\bfv_4.
       \hspace*{1cm}
       \nonumber
\end{eqnarray}
\QEDB
\end{lem}
Properties \eqref{eq:cpr} indicate that reactive collision terms assure the correct chemical exchange rates
for the considered chemical reaction \eqref{eq:cr}.


\subsection{Equilibrium state and $\cal H$-theorem}
\label{ssec:es}

In paper \cite{DMS-EJM-05}, the equilibrium solutions of the kinetic equations \eqref{eq:ke}
were studied in two steps, namely by considering first the mechanical equilibrium
associated to elastic collision operators
and then the chemical equilibrium associated to the reactive collisional operator.

The mechanical equilibrium is defined by distribution functions
{$f_i (\bfv , I )$,} $i=1,\ldots,4$,
such that
\begin{equation}
Q^e_{ij} (f_i,f_j) = 0,
\qquad i,j = 1,\ldots,4.
\label{eq:meq}
\end{equation}
Assuming that the mechanical equilibrium is reached,
the chemical equilibrium is then defined in paper~\cite{DMS-EJM-05} by distribution functions satisfying,
besides conditions \eqref{eq:meq}, the further condition
\begin{equation}
Q_i^{react} = 0, \qquad i=1,\ldots,4.
\label{eq:ceq}
\end{equation}
The following proposition summarizes the $\cal H$-theorems stated in paper \cite{DMS-EJM-05}
and characterizes the equilibrium states of the model.
\begin{prop}
(See Propositions $1$ and $2$, pages $231$--$233$ of $\cite{DMS-EJM-05}$)
\\
Let
$B_{ij}(\bfv,{\bfv}_{j},I,I_{j},R,r,\omega)$ and $B^{react}(\bfv,{\bfv}_j, I,I_j,R,r,\omega)$
be strictly positive almost everywhere
and let the distribution functions $f_i$ be non-negative for all $i=1,\ldots,4$
and such that the collisional operators $Q^e_{ij}$ and $Q_i^{react}$  are well defined.
\begin{enumerate}
\item Concerning the mechanical equilibrium, the following three properties are equivalent.
         \begin{enumerate}
         \item $Q_{ij}^e (f_i,f_j) = 0$
                 for all $i,j=1,\ldots,4$, \ $\bfv\in\RR^3$, \ $I\in\RR^+$;

                  \medskip

         \item $\dps \sum_{i,j=1}^4 \! {\int_{\RR^3}}  \! \int_0^{+\infty} \! Q_{ij}^e (f_i,f_j)
                           \log \! \left( \frac{f_i(\bfv,I)}{m_i^3} \right) \varphi_i(I) \, \dd I \dd\bfv = 0$;

                  \medskip

         \item There exist $n_i\geq0$, $i=1,\ldots,4$, $\bfu\in\RR^3$ and $T>0$ such that
                  $f_i(\bfv, I) =M( n_i , \bfu , T)$, where
                  \begin{equation}
                  \!\!
                  M( n_i , \bfu , T) \!=\! \frac{n_i(t,{\bfx})}{q_i(T(t,{\bfx}))} \!
                             \left( \! \frac{m_i}{2\pi k_BT(t,{\bfx})} \! \right )^{3/2} \!\!
                             \exp \! \left( \!\! -\frac{m_i\left |\bfv \!-\! {\bfu}(t,{\bfx}) \right |^2}{2k_BT(t,{\bfx})}
                             \!-\! \frac{I}{k_BT(t,{\bfx})} \! \right )
                  \label{eq:max}
                  \end{equation}
                  and
                  \begin{equation}
                  q_i(T(t,{\bfx}))  =  \int_0^{+\infty}\varphi_i(I) \, \exp \! \left( \! -\frac{I}{k_BT(t,{\bfx})} \! \right) \dd I,
                  \label{eq:qi}
                  \end{equation}

                  i.e. $f_i$ are Maxwellian distributions.
         \end{enumerate}

\item Assuming that the mechanical equilibrium is reached,
         that is, the distribution functions $f_i$ are given by \eqref{eq:max},
         the following three properties related to chemical equilibrium are equivalent
         \begin{enumerate}
         \item $Q_i^{react} = 0$
                  for all \ $i=1,\ldots,4$, \ $\bfv\in\RR^3$, \ $I\in\RR^+$;

                  \medskip

         \item $\dps \sum_{i=1}^4 \! \int_{\RR^3} \! \int_0^{+\infty} \! Q_i^{react}
                           \log \! \left( \frac{f_i(\bfv,I)}{m_i^3} \right) \varphi_i(I) \, \dd I \dd\bfv = 0$;

                  \medskip

         \item The following mass action law holds
                  \begin{equation}
                  \frac{n_1 n_2}{n_3 n_4} = \left( \frac{m_1 m_2}{m_3 m_4} \right)^{3/2}
                  \frac{q_1(T) q_2(T)}{q_3(T) q_4(T)} \; \exp\left( \frac{E}{k_BT} \right) .
                  \label{eq:mal}
                  \end{equation}
         \end{enumerate}
\end{enumerate}
\QEDB
\label{eq:ph}
\end{prop}
\section{Scaled equations and assumptions}
\label{sec:sea}
In this section we define the scaling regime for the kinetic equations
and introduce the assumptions to be considered at the kinetic level,
in order to derive the reaction-diffusion system of Maxwell-Stefan type as the hydrodynamic limit
equations of the considered kinetic model.

The evolution domain of the mixture is here represented by an open bounded domain $\Omega \subset \RR^3$,
with regular boundary.
\subsection{Scaling regime}
\label{ssec:seqs}
Let $\varepsilon$ be a scaling parameter  representing the mean free path or, equivalently, the Knudsen number,
with
$0 < \varepsilon \ll 1$. We scale the time and space variables as $(t,\bfx)\mapsto(\varepsilon^2 t , \varepsilon \bfx)$
and the distribution functions in the transformed variables are denoted by $f_i^\varepsilon$.

We consider a diffusive scaling regime for which
elastic collisions are dominant with respect to reactive collisions.
Accordingly, we start our analysis from the following scaled Cauchy problem for the distribution
functions $f_i^\varepsilon$,
\begin{align}
&  \varepsilon \frac{\partial f_i^\varepsilon}{\partial t} + \bfv \cdot \nabla_{\bfx} f_i^\varepsilon
       = \frac{1}{\varepsilon} \sum_{j=1}^4 Q_{ij}^e(f_i^\varepsilon,f_j^\varepsilon) + \varepsilon \, {Q_i}^{react} ,
       \quad i=1,\ldots,4 ,
       \label{eq:ske} \\[1.5mm]
& \hspace*{5.75cm}
       \mbox{for} \quad (t,\bfx,\bfv,I)\in \, ]0,+\infty[ \, \times\RR^3 \times \RR^3 \times [0,+\infty[ ,
        \nonumber \\[3mm]
&     f_i^\varepsilon(0 , \bfx , \bfv, I) = (f_i^{\rm in})^\varepsilon (\bfx , \bfv , I) ,
        \quad
        \mbox{for} \quad  (\bfx,\bfv,I)\in \RR^3 \times \RR^3 \times [0,+\infty[ .
        \label{eq:sic}
\end{align}
The properties of the collisional operators stated in Section \ref{sec:meqs},
with the obvious adjustments,
are still valid for the scaled operators.
\subsection{Assumptions}
\label{ssec:assump}
In view of the reaction-diffusion limit to be investigated in this paper,
we will consider the following assumptions on global mean velocity and temperature:
\begin{enumerate}
\item[{\it (a)}] The reactive mixture is a non-isothermal system,
               meaning that the temperature of the mixture is, in general,
               not constant in time and non-uniform in space.

\item[{\it (b)}] The bulk velocity of the mixture is small and goes to zero
               as the parameter $\varepsilon$ tends to zero.

\item[{\it (c)}] We assign as initial conditions $(f_i^{\rm in})^\varepsilon$
                     Maxwellian functions centered at the species mass velocity $\varepsilon{\bfu}_i^{\rm (in)}(\bfx)$
                     and species temperature $T^{\rm (in)}(\bfx) + \varepsilon T_i^{\rm (in)}(\bfx)$:
                     \begin{eqnarray}
                     (f_i^{\rm in})^\varepsilon (\bfx , \bfv , I)
                     \!\!&\!=\!&\!\! M_i^{\rm in}  \!=\!   
                     \frac{n^{\rm (in)}_i {(\bfx)}}{q_i\Big( \! T^{\rm (in)}(\bfx) + \varepsilon T_i^{\rm (in)}(\bfx) \! \Big)} \!
                     \left( \! \frac{m_i}{2\pi k_B \Big(T^{\rm (in)}(\bfx) + \varepsilon T_i^{\rm (in)}(\bfx)\Big)} \! \right )^{\!\!3/2} \!\!
        \nonumber
        \\[2mm]
        & & \times        
        \exp \! \left[ \! -\frac{m_i\left |\bfv \!-\! \varepsilon{\bfu}_i^{\rm (in)}(\bfx) \right |^2}{2k_B (T^{\rm (in)}(\bfx) + \varepsilon T_i^{\rm (in)}(\bfx))}\right]
        \exp \! \left[ \! - \frac{I}{k_B \Big(T^{\rm (in)}(\bfx) + \varepsilon T_i^{\rm (in)}(\bfx) \Big)} \right]   
        \hspace*{0.5cm}
\label{eq:infi}
\end{eqnarray}
for some
$n^{\rm (in)}_i ,\ T^{\rm (in)},\  T_i^{\rm (in)} \!\! : \, \Omega \rightarrow \RR^+$ and
${\bfu}_i^{\rm (in)}\!\! : \Omega \rightarrow \RR^3$.

\item[{\it (d)}] We assume that the time evolution of $f_i^\varepsilon$ established by equations \eqref{eq:ske}
               preserves for any time $t>0$ the initial Maxwellian structure of species distributions; more precisely, we consider that the distribution functions $f_i^\varepsilon$, at time $t>0$, are Maxwellians of the form
\begin{multline}
{f_i^\varepsilon(\bfv, I)} \!=\!
        \frac{n^\varepsilon_i(t,{\bfx})}{q_i\Big(T(t,{\bfx})+\varepsilon T_i^\varepsilon(t,{\bfx})\Big)}
        \Bigg( \frac{m_i}{2\pi k_B \Big(T(t,{\bfx})+\varepsilon T_i^\varepsilon(t,{\bfx})\Big)} \Bigg )^{\!\!3/2} \\
        \times 
        \exp \! \left[ \! -\frac{m_i\left |\bfv - \varepsilon{\bfu}_i^\varepsilon(t,{\bfx}) \right |^2}
        {2k_B \Big(T(t,{\bfx})+\varepsilon T_i^\varepsilon(t,{\bfx})\Big)} \right]
         \exp \left[ \! - \frac{I}{k_B \Big(T(t,{\bfx})+\varepsilon T_i^\varepsilon(t,{\bfx}) \Big)} \right]
\label{eq:maxhip}
\end{multline}
for some functions
$n^\varepsilon_i ,\ T_i^\varepsilon ,\ T \,: \,\RR^+ \times \Omega \rightarrow \RR^+$ and
${\bfu}_i^\varepsilon\, : \, \RR^+ \times \Omega \rightarrow \RR^3$.
\end{enumerate}
Assumptions {\it (c)} and {\it (d)} are consistent with the fact that the scaled equations \eqref{eq:ske} provide, as $\varepsilon\to0$,
$$
\sum_{j=1}^4 Q_{ij}^e(f_i^\varepsilon,f_j^\varepsilon) = O(\varepsilon)\,.
$$
Hence the distributions $f_i^\varepsilon$ should be, as $\varepsilon\to 0$, $O(\varepsilon)$ perturbations of collision equilibria of global elastic scattering operator, which are provided by local Maxwellians sharing the common zero mass velocity and a common temperature.
Therefore, for any time $t>0$, the species distributions may be considered to have the form
\begin{equation}
f_i^\varepsilon(\bfv, I) = M_i \Big( n_i(t,\bfx) , {\bfze} , T(t,\bfx) \Big)\Big(1 + O(\varepsilon) \Big),
\label{eq:dd}
\end{equation}
where $M_i$ is an elastic  {Maxwellian equilibrium state} as defined in \eqref{eq:max}.
Hence, up to the order $O(\varepsilon) $, we can suppose that
$$
f_i^\varepsilon(\bfv, I)=M_i \Big( n_i^\varepsilon(t,\bfx) , \varepsilon \bfu_i^\varepsilon(t,\bfx) , T(t,\bfx) + \varepsilon T_i^\varepsilon(t,\bfx) \Big)
$$
with
\begin{equation}
\lim_{\varepsilon\to0}
\Big( n_i^\varepsilon(t,\bfx) , \varepsilon \bfu_i^\varepsilon(t,\bfx) , T(t,\bfx) + \varepsilon T_i^\varepsilon(t,\bfx) \Big) =
\Big( n_i (t,\bfx), {\bfze} , T(t,\bfx)\Big),
\label{eq:lim}
\end{equation}
and
$$
\bfu_i^\varepsilon(t,\bfx) = O(1),\qquad T_i^\varepsilon(t,\bfx) = O(1).
$$
Assumption \eqref{eq:maxhip} is one of the simplest options to fulfill the property \eqref{eq:dd},
and it will be sufficient to obtain a system of Maxwell-Stefan type in the asymptotic limit $\varepsilon \to 0$.

\section{The limiting equations for the reactive mixture}
\label{sec:limin}
In this section, we derive the macroscopic equations in the hydrodynamic limit of the scaled kinetic equations
\eqref{eq:ske} with initial conditions defined by \eqref{eq:sic} and \eqref{eq:infi}.
As usual in kinetic theory \cite{ BGL91, BGL93, BD06, BGP-NA-17, HS-NA-17, HS-MMAS-17},
such equations are obtained by taking the appropriate moments of equations \eqref{eq:ske}
with respect to the velocity $\bfv$ and here also with respect to internal energy parameter~$I$.
In the present case, the balance equations obtained in the limit for $\varepsilon\to0$ constitute a reaction-diffusion system
of Maxwell-Stefan type for the reactive mixture of polyatomic gases.
The system is formed by the number density equation and momentum equation for each constituent,
as well as the balance equation for the mixture temperature.
Assumptions {\it (b)} and {\it (d)} considered in Subsection \ref{ssec:assump} play a crucial role
in the passage from the kinetic equations to the Maxwell-Stefan setting.


\subsection{Preliminaries}
\label{ssec:pre}

In what follows, instead of the unit vector $\omega\in \mathbb{S}^{2}$  parametrization used in \eqref{eq:KE_parametrization} and \eqref{eq:KE_parametrization1} to obtain the elastic and reactive post-collisional velocities,  we will use the  following unit vector $\sigma\in \mathbb{S}^{2}$ parametrization and show how to pass from one to the other.


\begin{prop}
Let
\begin{equation}
{\sigma} = {\bfy} - 2 (\omega\cdot{\bfy})\omega, \qquad
{\bfy}      = \frac{{\bfv} - {\bfv}_*}{\left |{\bfv}-{\bfv}_* \right |} ,
\label{eq:transformation}
\end{equation}
where ${\sigma}$ and $\omega$ are unit vectors in the sphere $\mathbb{S}^2$ {and $\bfv$, ${\bfv}_*$ denote the velocities of the colliding pair of particles.}
The Jacobian of the transformation from $\omega$ to $\sigma$ is given by
\begin{equation}
J_{{\omega} \rightarrow \sigma} = 4\left |\omega\cdot{\bfy}  \right |,  \qquad\forall \,{\bfy}\in\mathbb{R}^3.
\label{eq:Jacobian}
\end{equation}
\label{prop:sigma}
\end{prop}

\noindent
{\em Proof.}  See \cite{PThesis}, page 41.
\QEDB

\smallskip


\begin{rem} $\,$
\begin{enumerate}[label=(\roman*)]
\item Using the $\sigma$-parametrization given in \eqref{eq:transformation},
the post-collisional velocities \eqref{eq:elastic_bispecies_postcollisional_velocity} for a bi-species elastic collision between a particle of species~$i$ with ingoing parameters $(\bfv, I)$ and a particle of species~$j$ with parameters $(\bfv_j, I_j)$
can be rewritten as
\begin{equation}
\bfv_i'=\frac{m_i\bfv+m_j{\bfv}_j}{m_i+m_j}+\frac{m_j}{m_i+m_j}\sqrt{\frac{2R{\cal E}}{\mu_{ij}}} \, \sigma ,
\qquad
{\bfv}'_j=\frac{m_i{\bfv}+m_j{\bfv}_j}{m_i+m_j}-\frac{m_i}{m_i+m_j}\sqrt{\frac{2R{\cal E}}{\mu_{ij}}} \, \sigma .
\label{eq:v_iv_jprime}
\end{equation}
Therefore,
\begin{align}
\bfv_i'-\bfv
      & =\frac{m_j}{m_i+m_j}{\Vv}+\frac{m_j}{m_i+m_j}\sqrt{\frac{2R{\cal E}}{\mu_{ij}}}{\sigma},
      \label{eq:postv_minus_prev}
\end{align}
where
\begin{equation}
{\Vv}={\bfv}_j-\bfv.
\label{eq:relative_velocity}
\end{equation}

\item As a consequence of Proposition $\ref{prop:sigma}$,
         the bi-species elastic and reactive collision kernels can be written as, see $\cite{CV-02}$,
\begin{align}
& B_{ij}(\bfv,{\bfv}_j,I,I_j,R,r,\omega)\dd\omega= 2\left |\omega\cdot{\bfy} \right |
         B_{ij}(\bfv,{\bfv}_j,I,I_j,R,r,\sigma)\dd \sigma,
     \label{eq:transformation_elastic_kernel}\\[2mm]
& B^{react}(\bfv,{\bfv}_j,I,I_j,R,r,\omega) \dd\omega = 2\left |\omega\cdot{\bfy} \right |
        B^{react}(\bfv,{\bfv}_j,I,I_j,R,r,\sigma)\dd \sigma ,
\label{eq:transformation_reactive_kernel}
\end{align}
{where $\bfy = -\, \Vv/ |\Vv|$.}
\QEDB
\end{enumerate}
\end{rem}


\noindent
From now on we will sometimes skip for brevity the dependencies of distributions on velocity and internal energy, and we set $f^\varepsilon_i = f^\varepsilon_i(\bfv, I)$, $f^\varepsilon_j = f^\varepsilon_j (\bfv_j, I_j)$, $f^{\varepsilon'}_i = f^\varepsilon_i (\bfv_i', I_i')$, $f^{\varepsilon'}_j = f^\varepsilon_j(\bfv_j', I_j')$, and analogously for Maxwellian distributions.

\begin{prop}
\label{prop:BiRe}
Using the {the Taylor  expansion with respect to $\varepsilon$ of the distribution functions \eqref{eq:maxhip},} one can write
\begin{multline}
\mbox{\it(a)} \quad
{f}^{\varepsilon'}_i{f}^{\varepsilon'}_j-f^\varepsilon_if^\varepsilon_j
={M_i^\varepsilon M_j^\varepsilon\bigg[\varepsilon \Bigg( {\ba}_{ij}\cdot({\bfv}'_i- \bfv)\Bigg)
+\varepsilon\bigg(\frac{ T_i^\varepsilon}{k_B T^2}-\frac{T_j^\varepsilon}{k_B T^2}\bigg)} \\
\times \bigg( \frac{ m_i(v'_i)^2}{2}-\frac{ m_i v^2}{2}+I_i'-I\bigg)+O( \varepsilon^2)\bigg],
\label{eq:gain_minus_loss_bispecies}
\end{multline}
\begin{equation}
\begin{aligned}
\mbox{\it(b)} \quad
\left (\frac{m_i m_j}{m_k m_l} \right )^{\!3}
       & f^{\varepsilon'}_k f^{\varepsilon'}_l-f^\varepsilon_i f^\varepsilon_j
       \!=\! \left( \frac{m_i m_j}{m_k m_l} \right)^{\!3}
          M_k^{\varepsilon'}  M_l^{\varepsilon '}   -  M_i^\varepsilon M_j^\varepsilon
          \\
         & \!+\! \varepsilon \!\left[ \!\left(\!\frac{m_i m_j}{m_k m_l}\! \right)^{\!3}\!
          M_k^{\varepsilon'} M_l^{\varepsilon'}\!\bigg( \ba_k \!\cdot\! {\bfv}'_k
          +\frac{ m_k T_k^\varepsilon (v_k')^2}{2k_B T^2}+
          \frac{T_k^\varepsilon I_k'}{k_B T^2}
         -  \frac{T_k^\varepsilon}{T} \frac{q_k^*(T)}{k_BT q_k(T)}
         -  \frac{3T_k^\varepsilon}{2T}
          \right.
          \\
          & \qquad \left.
          + \ba_l \!\cdot\! {\bfv}'_l
          +\frac{ m_l T_l^\varepsilon (v_l')^2}{2k_B T^2}\!+\!
          \frac{T_l^\varepsilon I_l'}{k_B T^2}
          -  \frac{T_l^\varepsilon}{T} \frac{q_l^*(T)}{k_BT q_l(T)}
          -  \frac{3T_l^\varepsilon}{2T}
          \bigg) \right.
          \\
          & \qquad   \! - \!  \left.
          M_i^\varepsilon M_j^\varepsilon \bigg( \ba_i \!\cdot\! \bfv+
          \frac{ m_i T_i^\varepsilon v^2}{2k_B T^2}+\frac{T_i^\varepsilon I}{k_B T^2}
          -  \frac{T_i^\varepsilon}{T} \frac{q_i^*(T)}{k_BT q_i(T)}
          -  \frac{3T_i^\varepsilon}{2T}
          \right.
          \\
          & \qquad \left.+\ba_j \!\cdot\! {\bfv}_j+\frac{ m_j T_j^\varepsilon v_j^2}{2k_B T^2}+
          \frac{T_j^\varepsilon I_j}{k_B T^2}
          -  \frac{T_j^\varepsilon}{T} \frac{q_j^*(T)}{k_BT q_j(T)}
          -  \frac{3T_j^\varepsilon}{2T} \bigg)
          \!\right] \!
          + \! O(\varepsilon^2) ,
 \label{eq:postf-pref_reactive}
 \end{aligned}
 \end{equation}
where
\begin{equation}
{\ba}_{ij} = \frac{ m_i({\bfu}_i^\varepsilon-{\bfu}_j^\varepsilon)}{k_BT},\quad
{\ba}_i =  \frac{m_i{\bfu}^\varepsilon_i}{k_BT} ,\quad
{\ba}_j= \frac{m_j{\bfu}^\varepsilon_j}{k_BT} , \quad
{\ba}_k =  \frac{m_k{\bfu}^\varepsilon_k}{k_BT}, \quad
{\ba}_l =  \frac{m_l{\bfu}^\varepsilon_l}{k_BT} ,
\label{eq:ai}
\end{equation}
\begin{align}
q_i^*(T) &= \int_0^{+\infty} I \, \varphi_i(I) \, \exp \! \left( \! - \frac{I}{k_B T} \right )\dd I,
\label{eq:qistar}
\\
M_i^\varepsilon &:=
\frac{n_i^{\varepsilon}}{q_i(T)}
\left( \frac{m_i}{2\pi k_B T}  \right)^{3/2} \exp \left( -\frac{m_i v^2}{2k_BT}-\frac{I}{k_B T}  \right) ,
\label{eq:M}
\end{align}
and similarly for $M_j^\varepsilon$, $ M_k^{\varepsilon'}$, $M_l^{\varepsilon'}$.
\end{prop}
\noindent
{\em Proof.}  See the Appendix, part (A).
\QEDB

\medskip


\noindent 
Using $(a)$ of Proposition \ref{prop:BiRe} together with \eqref{eq:postv_minus_prev} {and \eqref{eq:transformation_elastic_kernel}},
we can write the elastic bi-species operator as
\begin{align}
Q_{ij}^e{(f_i,f_j)}   &\!=\! \varepsilon 
             \frac{m_j}{m_i \!+\! m_j}
             \int_{\mathbb{R}^3} \! \int_0^{+\infty} \!\! \int_0^1 \! \int_0^1 \! \int_{\mathbb{S}^2} \!
             (\ba_{ij} \!\cdot\! {\Vv}) \, M_i^\varepsilon M_j^\varepsilon \;
             {B_{ij}(\bfv,{\bfv}_j,I,I_j,R,r,\sigma)}
             \nonumber \\
             & \hspace*{5cm}
             \times \frac{2\cos\theta(1\!-\!R)}{V\varphi_i(I)}  \,\dd \sigma\,\dd r\, \dd R\, \dd I_j\, \dd{\bfv}_j
             \nonumber  \\
             &+\!\varepsilon   
             \frac{m_j}{m_i \!+\! m_j} \, {\ba_{ij}} \!\cdot\!
             \int_{\mathbb{R}^3} \! \int_0^{+\infty} \!\! \int_0^1 \! \int_0^1 \! \int_{\mathbb{S}^2} \!
             \sqrt{\frac{2R{{\cal E}}}{\mu_{ij}}} {\sigma}
             M_i^\varepsilon M_j^\varepsilon  \,
             {B_{ij}(\bfv,{\bfv}_j,I,I_j,R,r,\sigma)}
             \label{eq:Q_i_bispecies}
             \\
             & \hspace*{5cm}
             \times \frac{2\cos\theta(1\!-\!R)}{V\varphi_i(I)}\,\dd \sigma\,\dd r\, \dd R\, \dd I_j\, \dd{\bfv}_j
             \nonumber  \\
             &  {+\!\varepsilon
             \int_{\mathbb{R}^3} \! \int_0^{+\infty} \!\! \int_0^1 \! \int_0^1 \! \int_{\mathbb{S}^2} \!
             \bigg(\frac{ T_i^\varepsilon}{k_B T^2}-\frac{T_j^\varepsilon}{k_B T^2}\bigg)
             \bigg( \frac{ m_i(v'_i)^2}{2}-\frac{ m_i v^2}{2}+I_i'-I\bigg)M_i^\varepsilon M_j^\varepsilon \nonumber}\\
              & \hspace*{5cm}{\times{B_{ij}(\bfv,{\bfv}_j,I,I_j,R,r,\sigma)}\frac{2\cos\theta(1\!-\!R)}{V\varphi_i(I)}\,
              \dd \sigma\,\dd r\, \dd R\, \dd I_j\, \dd{\bfv}_j}
             \nonumber  \\[-2mm]
             & + O(\varepsilon^2) .
             \nonumber
\end{align}

\noindent
Similarly, using $(b)$ {of} Proposition \ref{prop:BiRe} {together with \eqref{eq:transformation_reactive_kernel}}, we can write the reactive operator as
\begin{multline}
Q_i^{react} =   \int_{\mathbb{R}^3}\int_0^{+\infty}\!\! \int_0^1\int_0^1\int_{\mathbb{S}^2}
                  \left[ \left(\frac{m_i m_j}{m_k m_l} \right)^{\!\!3} \!\!
                  M_k^{\varepsilon'} M_l^{\varepsilon '} - M_i^\varepsilon M_j^\varepsilon\right ] \!
                  B^{react}(\bfv,{\bfv}_j,I,I_j,R,r,\sigma)
                  \label{eq:Q_i_reactive}
                  \\
                  \times H_i\,
                  2\cos\theta(1-R)V^{-1}\varphi_i(I)^{-1} (m_im_j)^{-2}
                  \; \dd \sigma \dd r \dd R \dd I_j \dd{\bfv}_j
                  + O(\varepsilon) ,
\end{multline}
where
$\left |\omega\cdot{\bfy}  \right |= \cos\theta$, $V=\left | \bfv-{\bfv}_j \right |=\left \| {\Vv} \right \|$
{and for the forward and backward reactive operators, the {indices}
$(i,j,k,l)$ are such that  $(i,j,k,l) \in \left \{(1,2,3,4),(2,1,4,3) \right \} $ and $(i,j,k,l) \in \left \{(3,4,1,2),(4,3,2,1) \right \} $, respectively.}
The notation $M_i^\varepsilon , M_j^\varepsilon $ indicates Maxwellian distributions
defined as in \eqref{eq:M}.

\medskip

Also, we consider general collision kernels and split both the bi-species and reactive collision kernels
into the product of kinetic and angular collision kernels (see \cite{CV-02}) as given below
\begin{eqnarray}
& & B_{ij}(\bfv,{\bfv}_j,I,I_j,R,r,\sigma)
        = \underbrace{\left | \bfv-{\bfv}_j\right |^{\gamma} \Phi_{ij}(I,I_j,R,r)}_{\rm kinetic}\,
        \underbrace{b_{ij}(\cos\theta)}_{\rm angular} \, ,
        \label{eq:kernesplit}
        \\[2mm]
& & B^{react}(\bfv,{\bfv}_j,I,I_j,R,r,\sigma)
        = \underbrace{\left | \bfv-{\bfv}_j\right |^{\gamma} \Phi^{react}(I,I_j,R,r)}_{\rm kinetic}\,
        \underbrace{b^{react}(\cos\theta)}_{\rm angular} \, ,
         \label{eq:kernsplit}
\end{eqnarray}
where $\gamma$ is a parameter such that $\gamma \geq 1$.
This decomposition of the collision kernels is enough to obtain appropriate expressions for the integral production terms.

\medskip
Finally, we will use in the sequel the integral representation of both the gamma function
and the incomplete gamma function,
respectively defined as given below
\begin{align}
\int_0^{{+\infty}} x^n e^{-\alpha x^2}dx  &=  \frac{1}{2} \, \Gamma \left( \frac{n+1}{2} \right) \! \left(\frac{1}{\alpha } \right)^{\frac{n+1}{2}},
\label{eq:int_gamma_function}
\\[2mm]
\Gamma(\alpha,x)&=\int_x^{{+\infty}} t^{\alpha-1}e^{-t} dt.
\label{eq:int_incomplete_gamma_function}
\end{align}
\subsection{Moment of order zero}
\label{ssec:mom0}
From the scaled equations \eqref{eq:ske}, we first derive the evolution equation
for the number density of each constituent.
The relevant result is the following.
\begin{lem}
The balance equation for the number density of each constituent in the reactive mixture
can be written in the form
\begin{equation}
\frac{\partial n_i }{\partial t} + \nabla_{\bfx}\cdot(n_i {\bfu}_i) = -\lambda_i A,
\quad i=1,2,3,4,
\label{eq:Conservation_of_mass_final}
\end{equation}
where
\begin{equation*}
n_i = \lim_{\varepsilon\rightarrow 0}n_i^\varepsilon, \qquad
{\bfu}_i = \lim_{\varepsilon\rightarrow 0}{\bfu}_i^\varepsilon, \qquad
\end{equation*}
and
\begin{equation}
\lambda_1=\lambda_2=1, \qquad
\lambda_3=\lambda_4= - 1
\label{eq:lamb}
\end{equation}
are stoichiometric coefficients.
Moreover, $A$ is the production term defining the chemical rate
of the reactants of the backward 
chemical reaction.
More specifically, $A$ is the formal limit for $\varepsilon\to0$ of
\begin{eqnarray}
A^\varepsilon \!\!&\!\!=\!\!&\!\!
                 \left[\!
                 \bigg(\frac{m_3m_4}{m_1m_2}\bigg)^{\!3/2}
                 \frac{n_1^\varepsilon n_2^\varepsilon}
                 {q_1(T) q_2(T)}
                 \exp \left( - \frac{E}{k_B T} \! \right)
                 -
                 \frac{n_3^\varepsilon n_4^\varepsilon}
                 {q_3(T) q_4(T)} \! \right]
                  \label{eq:Aetauu}
                  \\[2mm]
        & &   \times
                 \frac4{\sqrt\pi(m_3 m_4)^2}  \left( \frac{2 k_B T}{\mu_{34}} \right)^{\!\!(\gamma-1)/2}
                 \Gamma \! \left( \frac{\gamma + 2}{2} \right) \!
                 \int_0^{+\infty} \!\! \int_0^{+\infty} \!\!
                 \int_0^1 \! \int_0^1 \! \int_{\mathbb{S}^2} \!
                 \exp \! \left( -\frac{I\!+\!I_4}{k_B T} \!\right)
                 \nonumber
               \\[3mm]
        & &   \times {\Phi^{react}(I , I_4 , R , r)}
                 b^{react}(\cos\theta) {\cos \theta} \, (1-R) \,
                 \dd \sigma \dd r \dd R \dd I_4 \dd I + O(\varepsilon).
                 \nonumber
\end{eqnarray}
\label{lm:mom0}
\end{lem}

\medskip

\noindent
{\it Proof.}
Multiplying equations \eqref{eq:ske} by $\varphi_i(I)$ and integrating with respect to
$\bfv \!\in\! \mathbb{R}^3$ and $I \!\in\! \mathbb{R}_+$,
we obtain
\begin{eqnarray}
\lefteqn{
\varepsilon\frac{\partial }{\partial t}\underbrace{\left (\! \int_{\mathbb{R}^3}\int_0^{+\infty}f^\varepsilon_i \varphi_i(I)\dd I\,\dd\bfv\right )}_{n^\varepsilon_i}
            + \; {\nabla_{\bfx}\cdot}\underbrace{\left (\! \int_{\mathbb{R}^3}\int_0^{+\infty} \bfv
                               f_i^\varepsilon \varphi_i(I)\dd I \dd\bfv\right)}_{\varepsilon n^\varepsilon_i {\bfu}^\varepsilon_i} }
                               \hspace*{1.5cm}
            \label{eq:conservation_of_mass} \\[1mm]
   \!\!&\!\!=\!\!&\!\!  \frac{1}{\varepsilon} \underbrace{\int_{\mathbb{R}^3}\!\int_0^{+\infty}  \!\!
                    {Q_{ii}^e} \varphi_i(I)\dd I \dd\bfv}_{0} +
            \frac{1}{\varepsilon}{\sum_{\substack{j=1\\ j\neq i}}^{4}} \underbrace{\int_{\mathbb{R}^3}\int_0^{+\infty} \!\!
                    {Q_{ij}^e} \varphi_i(I)\dd I \dd\bfv}_{0}  +
            \varepsilon \underbrace{\int_{\mathbb{R}^3} \int_0^{+\infty} \!\!
                    Q_{i}^{react} \varphi_i(I) \dd I \dd\bfv}_{-\lambda_i A^\varepsilon}
                    \nonumber
\end{eqnarray}
where the term $A^\varepsilon$ is given by \eqref{eq:Aetauu}.
Moreover, we have used properties \eqref{eq:weakform_mono} and \eqref{eq:weakform_bispecies_collision_integral1},
respectively, with $\psi({\bfv}, I) \!=\!1$,
for the vanishing of the first two terms on the right hand side of equation~\eqref{eq:conservation_of_mass}.
To obtain the last term on the right hand side of  equation~\eqref{eq:conservation_of_mass},
we have used the more convenient form of the reactive operator $Q_{i}^{react}$, given in \eqref{eq:Q_i_reactive}.
Then we have used the decomposition of the reactive kernel given in \eqref{eq:kernsplit}
and the conservation of total energy for reactive collisions
given by condition\ \eqref{eq:conservation_eng_chemical}.
See the Appendix, {part{\it (B)}.}

\medskip

Finally, dividing both sides of equation~\eqref{eq:conservation_of_mass} by $\varepsilon$ and
taking the limit as $\varepsilon\rightarrow 0$,
we obtain the evolution equation \eqref{eq:Conservation_of_mass_final} for the constituent number density.
\QEDB


\subsection{Moment of order one}
\label{ssec:mom1}

From the scaled equations \eqref{eq:ske}, we then derive
the evolution equation for the momentum of each species,
as stated in the following result.
\begin{lem}
(i)
The balance equation for the momentum of each constituent in the reactive mixture
can be written in the form
\begin{equation}
{\nabla_{\bfx} (n_i k_B T) = m_i\bfB_i} ,  \quad i=1,2,3,4, \quad
n_i = \lim_{\varepsilon\rightarrow 0} n_i^\varepsilon,  \quad
{T = \lim_{\varepsilon\rightarrow 0} \big( T+\varepsilon T^\varepsilon_i \big)} ,
\label{eq:momentum_balance_final}
\end{equation}
where  
$\bfB_i$ is the production term associated to the momentum balance of the species,
namely the formal limit for $\varepsilon\to0$ of
\begin{equation}
\begin{aligned}
\bfB_i^\varepsilon &\!\!=\!\!  \sum_{\substack{j=1\\ j\neq i}}^{4}
           \frac{m_j}{m_i+m_j}
           \int_{\mathbb{R}^3} \!\! \int_{\mathbb{R}^3} \!\! \int_0^{+\infty} \!\!\! \int_0^{+\infty} \!\!\!
           \int_0^1 \!\! \int_0^1 \!\! \int_{\mathbb{S}^2} \!
           {(\ba_{ij}\!\cdot\! {\Vv})\bfv} M_i^\varepsilon M_j^\varepsilon
           V^{\gamma-1} \Phi_{ij}(I,I_j,R,r)
           \label{eq:B_i*}  \\[2mm]
           & \hspace*{3cm} \times 2 \, b_{ij}(\cos\theta) \cos\theta(1-R) \,
           \dd \sigma \dd r \dd R \dd I_j \dd I \dd{\bfv}_j \dd\bfv
            \\[2mm]
          &  +\!  \sum_{\substack{j=1\\ j\neq i}}^{4}
          \frac{m_j}{m_i+m_j}
          \int_{\mathbb{R}^3} \!\! \int_{\mathbb{R}^3} \!\! \int_0^{+\infty} \!\!\! \int_0^{+\infty} \!\!\!
          \int_0^1 \!\! \int_0^1 \!\! \int_{\mathbb{S}^2} \! (\ba_{ij} \! \cdot \! \sigma) \bfv
          \sqrt{\frac{2R {\cal E}}{\mu_{ij}}}
          M_i^\varepsilon M_j^\varepsilon V^{\gamma-1} \Phi_{ij}(I,I_j,R,r)
           \\[2mm]
           & \hspace*{3cm} \times 2 \, b_{ij} (\cos\theta) \cos\theta(1-R) \,
          \dd \sigma \dd r \dd R \dd I_j \dd I \dd{\bfv}_j \dd\bfv
          \\[2mm]
          & {+\! \sum_{\substack{j=1\\ j\neq i}}^{4}
             {\int _{{\mathbb{R}}^3}\!\int_{\mathbb{R}^3} \!\int_0^{+\infty}\!\!\! \int_0^{+\infty} \!\! \int_0^1 \! \int_0^1 \! \int_{\mathbb{S}^2}} \!
             \bfv\bigg(\!\frac{ T_i^\varepsilon}{k_B T^2}\!-\!\frac{T_j^\varepsilon}{k_B T^2}\!\bigg)\bigg( \frac{ m_i(v'_i)^2}{2}-\frac{ m_i v^2}{2}+I_i'-I\bigg)M_i^\varepsilon M_j^\varepsilon}
             \\[2mm]
        & \hspace*{3cm}{\times V^{\gamma-1} \Phi_{ij}(I,I_j,R,r) 2 \, b_{ij} (\cos\theta) \cos\theta(1-R) \,
          \dd \sigma \dd r \dd R \dd I_j \dd I \dd{\bfv}_j \dd\bfv}
             \\
           &+ O(\varepsilon) .
          \end{aligned}
\end{equation}

\noindent
(ii)
Moreover, if the angular collision kernel of the elastic bi-species kernel
is an odd function of $\cos\theta$,
then the production term $\bfB_i$ is the formal limit for $\varepsilon\to0$ of
\begin{equation}
\begin{aligned}
\bfB_i^\varepsilon\!\!=\! - & \! \sum_{\substack{j=1\\ j\neq i}}^{4} \! {\left(\frac{m_j}{m_i\!+\!m_j}\right)^2\!\!
      \bigg(\frac{m_i}{2\pi k_B T}\bigg)^{3/2}\!\!\bigg(\frac{m_j}{2\pi k_B T}\bigg)^{3/2}}
      \frac{n_j^\varepsilon{\bf{J}}_i^\varepsilon-n_i^\varepsilon{\bf{J}}_j^\varepsilon}{q_i(T)q_j(T)} \!
      \\[-2mm]
      & \times {\frac{m_i}{k_B T}}\left( \frac{2\pi k_B T}{m_i+m_j} \right)^{\!3/2}
      \frac{2\pi}{3}\Gamma \! \left( \! \frac{\gamma+4}{2} \! \right) \!
      \left( \! \frac{2k_B T}{\mu_{ij}} \! \right)^{\!\! \frac{\gamma+4}{2}}
      \int_0^{+\infty} \!\! \int_0^{+\infty} \!\! \int_0^1 \! \int_0^1 \! \int_{\mathbb{S}^2} \!\!
      \exp \! \left( \! -\frac{I \!+\!I_j}{k_B T} \! \right)
      \\[2mm]
       & \times \Phi_{ij}(I,I_j,R,r)
       2\cos\theta (1 \!-\! R) b_{ij}(\cos\theta) \, \dd \sigma \dd r \dd R \dd I_j \dd I \\
       &+ O(\varepsilon),
        \label{eq:Bbar_final_1}
\end{aligned}
\end{equation}

\noindent
where ${\bf{J}}_i^\varepsilon$ denotes the molar diffusive {flux} of the species
given by
\begin{equation}
{\bf{J}}_i^\varepsilon=n_i^\varepsilon({\bfu}_i^{\varepsilon}-{\bfu}^{\varepsilon}).
\label{eq:Ji}
\end{equation}
\label{lm:mom1}
\end{lem}
{\it Proof.}
(i)
First we multiply equations \eqref{eq:ske} by $\bfv\varphi_i(I)$ and integrate
with respect to $\bfv \!\in\! \mathbb{R}^3$ and $I\in \mathbb{R}_+$ to obtain
\begin{eqnarray}
\lefteqn{\varepsilon\frac{\partial }{\partial t} \underbrace{\left \{ \int_{\mathbb{R}^3} \!\! \int_0^{+\infty}\bfv
         f_i^\varepsilon \varphi_i(I) \dd I \dd\bfv \right \}}_{\varepsilon n_i^\varepsilon {\bfu}_i^\varepsilon} +
         {\nabla_{\bfx}\cdot}\underbrace{\left \{ \int_{\mathbb{R}^3} \!\! \int_0^{+\infty}\bfv\otimes \bfv
         f_i^\varepsilon \varphi_i(I) \dd I \dd\bfv \right \}}_{\frac{k_B{(T + \varepsilon T_i^\varepsilon)}}{m_i}
         n_i^\varepsilon +
         \varepsilon^2\left ( n_i^\varepsilon{\bfu}_i^\varepsilon \otimes {\bfu}_i^\varepsilon \right )} }
         \label{eq:momentum_balance}\\
         \!\!&\!\!=\!\!&\!\!
         \frac{1}{\varepsilon}\underbrace{\int_{\mathbb{R}^3} \!\!  \int_0^{+\infty}\bfv {Q_{ii}^e} \varphi_i(I) \dd I \dd\bfv}_{\small{\bf0}}
         + \underbrace{\frac{1}{\varepsilon}{\sum_{\substack{j=1\\ j\neq i}}^{4}} \int_{\mathbb{R}^3} \!\!  \int_0^{+\infty}\bfv {Q_{ij}^e} \varphi_i(I) \dd I \dd\bfv}_{ \bfB_i^\varepsilon}
         + \underbrace{\varepsilon \int_{\mathbb{R}^3} \!\!  \int_0^{+\infty}\bfv Q_{i}^{react} \varphi_i(I) \dd I \dd\bfv}_{ \bfC_i^\varepsilon}\, ,
\nonumber
\end{eqnarray}
where property \eqref{eq:weakform_mono} was used with $\psi(\bfv,I)=\bfv$
to obtain the vanishing of the first term on the right hand side.
Moreover, concerning the term $\bfB_i^\varepsilon$ on the right-hand side of \eqref{eq:momentum_balance}
we have used the
elastic bi-species operator $Q_{ij}^e$
given in \eqref{eq:Q_i_bispecies} and the decomposition of the elastic bi-species kernel indicated in \eqref{eq:kernesplit}
to obtain the expression for $\bfB_i^\varepsilon$ given in \eqref{eq:B_i*}.
Upon specifying $\Phi_{ij}(I,I_j,R,r)$,
the integrals with respect to $I,I_j,R,r,\sigma$ can be  evaluated explicitly
and the six fold integrals with respect to $\bfv, {\bfv}_j$ can be transformed to the center of mass velocity and relative velocity  and the resulting integrals can be evaluated using the integral representation of gamma function given in \eqref{eq:int_gamma_function}.

\medskip

\noindent
Then, concerning the term $\bfC_i^\varepsilon$ on the right-hand side of \eqref{eq:momentum_balance},
we have used the form of the reactive operator $Q_i^{react}$ given in \eqref{eq:Q_i_reactive}
and the decomposition \eqref{eq:kernsplit} for
the reactive kernel $B^{react}(\bfv,{\bfv}_j,I,I_j,R,r,\sigma)$
 to obtain
\begin{align}
\bfC_i^\varepsilon &\!=\!  \varepsilon \int_{\mathbb{R}^3}\int_{\mathbb{R}^3}
\int_0^{+\infty}\!\!\!\int_0^{+\infty}\!\!\! \int_0^1\int_0^1\int_{\mathbb{S}^2}
          \! \bfv \! \left [ \! \left(\frac{m_i m_j}{m_k m_l} \right)^{\!3} \!\!
          M_k^{\varepsilon'} M_l^{\varepsilon '} \!-\! M_i^\varepsilon M_j^\varepsilon\right ]
          V^{\gamma-1}\Phi^{react}(I, I_j,R,r)
          \label{eq:Ci1}\\[2mm]
          & \hspace*{2cm} \times 2 \, b^{react}(\cos\theta)H_i(m_im_j)^{-2} \cos\theta(1-R) \;
          \dd \sigma \dd r \dd R \dd I_j\dd I \dd{\bfv}_j\dd\bfv
         + O({\varepsilon^2}).
         \nonumber
\end{align}
Again, upon specifying $\Phi^{react}(I,I_j,R,r)$,
the integrals on the right hand side of equation~\eqref{eq:Ci1} can be computed.
Specifically, the integrals with respect to $I,I_j,R,r,\sigma$ can be evaluated explicitly
and the six fold integrals with respect to $\bfv, {\bfv}_j$ can be transformed to the center of mass velocity and relative velocity and the resulting integrals can either be evaluated using
the integral representation of gamma function given in \eqref{eq:int_gamma_function} or be represented by the incomplete gamma function given in \eqref{eq:int_incomplete_gamma_function}.

\medskip


\noindent
Finally we take the limit as $\varepsilon\rightarrow 0$ in equation~\eqref{eq:momentum_balance}
and obtain the balance equation \eqref{eq:momentum_balance_final} for the momentum of each species,
in which the production term $\bfB_i$  is the formal limit for $\varepsilon\to0$ of the term
$\bfB_i^\varepsilon$ given in~\eqref{eq:B_i*}.

\bigskip


\noindent
(ii)
Using the assumption that the angular collision kernel is an odd function of $\cos\theta$,
one obtains that the second term on the right hand side of equation~\eqref{eq:B_i*} vanishes.
See the {Appendix,} {part {\it (C1)}.}
Transforming the six fold integral over $\bfv$ and $\bfv_j$ in the first and third terms on the right hand side of \eqref{eq:B_i*} to the center of mass velocity and relative velocity and evaluating the resulting integrals using \eqref{eq:int_gamma_function},
we obtain that the third term on the right hand side of \eqref{eq:B_i*} vanishes.
See the Appendix, part {\it {(C3)}}.
Moreover,
the six fold integral over $\bfv$ and $\bfv_j$ in the first term on the the right hand side of the same equation reduces to
\begin{multline}
\int_{\mathbb{R}^3} \int_{\mathbb{R}^3} (\ba_{ij}\cdot {\bf{V}}) \bfv  M_i^\varepsilon M_j^\varepsilon V^{\gamma-1} \dd{\bfv}_j \dd\bfv
       =  - \frac{m_j}{m_i + m_j} \frac{n_i^\varepsilon n_j^\varepsilon}{q_iq_j}\frac{(m_im_j)^{3/2}}{(2\pi k_BT)^3} \;
       \exp \! \left( \! -\frac{I+I_j}{k_B T} \! \right ) \frac{m_i}{k_B T} ({\bfu}_i^{\varepsilon}-{\bfu}_j^{\varepsilon}) \\
       \times\left ( \frac{2\pi k_B T}{{m_i + m_j}} \right )^{3/2}\frac{2\pi}{3}\Gamma
       \left ( \frac{\gamma+4}{2} \right )\left ( \frac{2k_B T}{\mu_{ij}} \right )^{(\gamma+4)/2} .
\label{eq:int_vivj_Bbar}
\end{multline}
See the Appendix part {{\it (C2)}}.
Now, substituting \eqref{eq:int_vivj_Bbar} into the first addend in \eqref{eq:B_i*},
we then obtain the desired result.
\QEDB


\medskip

\begin{rem}
\label{rem:B_i_eta_reduction}
Observe that the production term $\bfB_i^\varepsilon$ given in equation~\eqref{eq:Bbar_final_1}
has the Maxwell-Stefan structure $\cite{BGP-NA-17, HS-MMAS-17}$.
The remaining integrals involved in equation~\eqref{eq:Bbar_final_1}, when evaluated,
will only contribute to the definition of the diffusion coefficients.
\QEDB
\end{rem}


\subsection{Conservation of energy}
\label{ssec:ce}

From the scaled equations \eqref{eq:ske}, we finally derive the evolution equation for the energy of the mixture
as stated in the following lemma.
\begin{lem}
The balance equation for the energy of the reactive mixture
can be written in the form
\begin{equation}
\frac{\partial }{\partial t} \left[ \sum_{i=1}^{4} \left( \frac{3}{2 } n_i k_BT  +  n_i \frac{q_i^*(T)}{q_i(T)} \right) \right]
      +  {\nabla_{\bfx}\cdot} \left[  \sum_{i=1}^{4} \left( \frac{5}{2} n_i k_BT \, {\bfu}_i
      +  n_i {\bfu}_i \frac{q_i^*(T)}{q_i(T)} \right) \right]
      = - E A ,
      \label{eq:moment_order_twoA}
\end{equation}
where
$$
n_i = \lim_{\varepsilon\rightarrow 0} n_i^\varepsilon ,  \qquad \qquad 
{\bfu}_i = \lim_{\varepsilon\rightarrow 0} {\bfu}_i^\varepsilon ,\qquad \qquad
{T = \lim_{\varepsilon\rightarrow 0} \big( T+\varepsilon T_i^\varepsilon\big)},
$$
 $q_i^*(T)$ has been defined in~\eqref{eq:qistar}, and $A$ is the production term defined as in Subsection $\ref{ssec:mom0}$.
\label{lm:mom2}
\end{lem}
{\it Proof.}
Multiplying Eqs.~\eqref{eq:ske} by $\left [ m_i v^2/{2}+I \right ]\varphi_i(I)$,
integrating with respect to $\bfv \! \in \! \mathbb{R}^3$ and $I \! \in \! \mathbb{R}_+$,
and then summing over the species $i=1,2,3,4$,
we obtain
\begin{eqnarray}
\lefteqn{ \varepsilon \frac{\partial }{\partial t}
       \sum_{i=1}^{4} \underbrace{\left\{
       \int_{\mathbb{R}^3} \int_0^{+\infty} \left[ \frac{ m_i v^2}{2} + I \right] \varphi_i(I) f_i^\varepsilon
       \dd I \dd\bfv \right\}}_{\frac{3}{2 }n_i^{\varepsilon}k_B {(T + \varepsilon T_i^\varepsilon)} +\frac{1}{2}\varepsilon^2 m_i n_i^{\varepsilon}({\bfu}_i^\varepsilon)^2
       + n_i^\varepsilon
       \frac{q_i^*({T + \varepsilon T_i^\varepsilon})}{q_i({T + \varepsilon T_i^\varepsilon})}}   }
       \hspace*{1cm}
       \label{eq:moment_order_two} \\[2mm]
& & \hspace*{1cm} + {\nabla_{\bfx}\cdot}
       \sum_{i=1}^{4}\underbrace{\left\{ \int_{\mathbb{R}^3} \int_0^{+\infty}
       \left[ \frac{ m_i v^2}{2}+I \right ]\varphi_i(I) \bfv f_i^\varepsilon
       \dd I \dd\bfv\right \}}_{\varepsilon\left \{ \frac{5}{2} n_i^\varepsilon k_B {(T + \varepsilon T_i^\varepsilon)}
       + \varepsilon^2 \left ( \frac{1}{2}m_in_i^\varepsilon({\bfu}_i^\varepsilon)^2 \right) \right \} {\bfu}_i^\varepsilon
       + \varepsilon  n_i^\varepsilon{\bfu}_i^\varepsilon
       \frac{q_i^*({T + \varepsilon T_i^\varepsilon})}{q_i({T + \varepsilon T_i^\varepsilon})}}
       = \;\;{\cal D}^\varepsilon + {\cal E}^\varepsilon + {\cal F}^\varepsilon,
       \nonumber
\end{eqnarray}
where ${\cal D}^\varepsilon$, ${\cal E}^\varepsilon$, ${\cal F}^\varepsilon$ are the production terms defined by
\begin{align*}
&  {\cal D}^\varepsilon = \frac{1}{\varepsilon}\sum_{i=1}^{4}\int_{\mathbb{R}^3}\int_0^{+\infty}
     \left[ \frac{ m_i v^2}{2}+I \right ]\varphi_i(I) {Q_{ii}^e}  \dd I \dd\bfv ,
     \\
&   {\cal E}^\varepsilon = \frac{1}{\varepsilon}\sum_{i=1}^{4}
{\sum_{\substack{j=1\\ j\neq i}}^{4}}
\int_{\mathbb{R}^3}\int_0^{+\infty}
     \left[ \frac{ m_i v^2}{2}+I \right ]\varphi_i(I) {Q_{ij}^e} \varphi_i(I)
     \dd I \dd\bfv ,
     \\
&  {\cal F}^\varepsilon = \varepsilon\sum_{i=1}^{4}\int_{\mathbb{R}^3}\int_0^{+\infty}
    \left[ \frac{ m_i v^2}{2}+I \right]
    \varphi_i(I) Q_{i}^{react}
    \dd I \dd\bfv .
\end{align*}
Using properties \eqref{eq:weakform_mono}, \eqref{eq:weakform_bispecies_collision_integral2} and \eqref{eq:weakform_reactive_collision_integral},
with $\psi(\bfv,I) \!=\! \left [  m_i v^2/2+I \right ]$,
for the {production} terms
${\cal D}^\varepsilon$, ${\cal E}^\varepsilon$ and ${\cal F}^\varepsilon$, respectively,
we obtain that
\begin{equation}
{\cal D}^\varepsilon=0, \qquad {\cal E}^\varepsilon=0 \quad{\text{and}}\quad {\cal F}^\varepsilon
   = - \varepsilon \, E A^\varepsilon ,
\label{eq:DEF}
\end{equation}
with $A^\varepsilon$ being defined by expression \eqref{eq:Aetauu}.
Substituting \eqref{eq:DEF} into  \eqref{eq:moment_order_two} and dividing both sides by $\varepsilon$, we obtain
\begin{eqnarray}
\lefteqn{  {\frac{\partial }{\partial t}} \left[ {\sum_{i=1}^4} \,
      \left( \frac{3}{2 }n_i^{\varepsilon}k_B {(T + \varepsilon T_i^\varepsilon)} +
      \frac{1}{2}\varepsilon^2 m_i n_i^{\varepsilon}({\bfu}_i^\varepsilon)^2
      + n_i^\varepsilon\frac{q_i^*({T + \varepsilon T_i^\varepsilon})}
      {q_i ({T + \varepsilon T_i^\varepsilon})} \right) \right]   }
      \label{eq:moment_order_twoB}\\
& & + {\nabla_{\bfx}\cdot} \left\{  {\sum_{i=1}^4}
       \left[ \left( \frac{5}{2} n_i^\varepsilon k_B {(T + \varepsilon T_i^\varepsilon)}
       +  \frac{1}{2} \varepsilon^2 m_i n_i^\varepsilon({\bfu}_i^\varepsilon)^2 \right) {\bfu}_i^\varepsilon
       +  n_i^\varepsilon{\bfu}_i^\varepsilon
       \frac{q_i^*({T + \varepsilon T_i^\varepsilon})}
       {q_i({T + \varepsilon T_i^\varepsilon})} \right] \right\} = - E A^\varepsilon .
       \nonumber
\end{eqnarray}
Taking the limit as $\varepsilon\to0$ in equation~\eqref{eq:moment_order_twoB},
we obtain the balance equation for the energy of the mixture in the form of
equation~\eqref{eq:moment_order_twoA}.
\QEDB


\subsection{Limit equations}
\label{ssec:limEqs}

In this section, we summarize the results of Section \ref{sec:limin}
and write the limit equations
when the angular collision kernels of both the elastic bi-species and the reactive operators are  odd functions of $\cos\theta$.
By recalling that, due to assumption $b)$, in the limit as $\varepsilon\to0$, we have $\bfJ_i \!=\! n_i\bfu_i$, for all {$i=1,2,3,4$,}
and putting together
the balance equations \eqref{eq:Conservation_of_mass_final}
and \eqref{eq:momentum_balance_final}
for the number densities and momentum of the constituents
as well as
the balance equation \eqref{eq:moment_order_twoA}
for the total energy of the mixture,
we obtain the following macroscopic system of reaction-diffusion equations,
\begin{eqnarray}
& & \frac{\partial n_i }{\partial t} + {\nabla_{\bfx}\cdot} \bfJ_i = -\lambda_i A,
       \quad  i=1,2,3,4,
       \nonumber
       \\[1mm]
& & \nabla_{\bfx}\Big(n_i k_B T \Big) =  -\sum_{\substack{j=1\\ j\neq i}}^{4}
       \frac{n_j{\bf{J}}_i-n_i{\bf{J}}_j}{D_{ij}},
       \quad  i=1,2,3,4,
       \label{eq:MS}
       \\[1mm]
& & \frac{\partial }{\partial t} \left[ \sum_{i=1}^{4} n_i  \left(\frac{3}{2 } k_B T  +  \frac{q_i^*(T)}{q_i(T)} \right) \right]
      +  {\nabla_{\bfx}\cdot} \left[ \sum_{i=1}^{4}\left( \frac{5}{2} k_BT
      + \frac{q_i^*(T)}{q_i(T)} \right) \! {\bfJ}_i  \right]= -E A,
      \nonumber
\end{eqnarray}
where $A$ is the production term derived in Subsection \ref{ssec:mom0}
as the formal limit for $\varepsilon \to 0$ of the
the approximate production term $A^\varepsilon$ given in equation~\eqref{eq:Aetauu}.
Moreover, $D_{ij}$, for {$i=1,2,3,4$,} and $i\not=j$, are diffusion coefficients
that can be recovered from the production term $\bfB_i$ obtained in Subsection \ref{ssec:mom1}
as the formal limit for $\varepsilon \to 0$ of $\bfB_i^\varepsilon$.
The production term $A$ is given by
\begin{eqnarray}
A &=&   \left[\left(\frac{m_3 m_4}{m_1 m_2}\right)^{\!\!3/2}
                 \frac{n_1 n_2}{q_1(T(t,\bfx)) q_2(T(t,\bfx))}
                 \exp \left( \! - \frac{E}{k_B T(t,\bfx)} \! \right) \!-\!
                 \frac{n_3 n_4 }{q_3(T(t,\bfx)) q_4 (T(t,\bfx))}  \right]
                 \nonumber
                 \\[2mm]
        & &   \times {\frac4{\sqrt{\pi}(m_3 m_4)^{2}} \, \Gamma\left( \frac{\gamma+2}{2}\right)}
                 \left( \frac{2k_BT}{\mu_{34}} \right)^{(\gamma-1)/2}
                 \int_{0}^{+\infty} \int_{0}^{+\infty}
                 \int_0^1 \! \int_0^1 \! \int_{\mathbb{S}^2} \!
                 \exp \! \left( -\frac{I\!+\!I_4}{k_B T} \!\right)
                                   \label{eq:Aend}
                                  \\[2mm]
        & &   \times
	        {\Phi^{react}(I , I_4 , R , r)}
                 b^{react}(\cos\theta) \, \cos \theta \, (1-R) \,
                 \dd \sigma \dd r \dd R \dd I_4 \dd I ,
                 \nonumber
\end{eqnarray}
and the diffusion coefficients $D_{ij}$ are defined as
\begin{align}
\frac1{D_{ij}} &=\frac{4}{3 \sqrt{\pi}}\,\mu_{ij}^{({3}-\gamma)/2}  (2k_B T(t,\bfx))^{(\gamma{-}1)/2}
      \frac{1}{q_i(T(t,\bfx))q_j(T(t,\bfx))}
      \Gamma \! \left( \! \frac{\gamma+4}{2} \! \right)
       \label{eq:Dij}
      \\[3mm]
      & \times \int_0^{+\infty} \! \int_0^{+\infty} \! \int_0^1 \!\! \int_0^1 \!\! \int_{\mathbb{S}^2} \!\!
      \exp \! \left( \! -\frac{I \!+\!I_j}{k_B T(t,\bfx)} \! \right) \! \Phi_{ij}(I,I_j,R,r)
      \cos\theta(1 \!-\! R) b_{ij}(\cos\theta) \, \dd \sigma \dd r \dd R \dd I_j \dd I \, .
      \nonumber
\end{align}

\medskip

\noindent
The equations in the second row of \eqref{eq:MS} are the Maxwell-Stefan equations
for the reactive mixture of polyatomic gases considered in this paper.
Upon specifying the kinetic and the angular kernels,
the integrals appearing in equations \eqref{eq:Aend} and \eqref{eq:Dij}
can be evaluated explicitly and detailed expressions for the reactive production terms $A$, $EA$
and diffusion coefficients $D_{ij}$ can be obtained.


\section{Conclusion}
\label{sec:concl}

In this article we have derived a set of reaction diffusion equations of Maxwell-Stefan type
for describing a chemically reactive gaseous mixture composed of polyatomic species,
which takes into account
both the presence of internal energy degrees as well as the chemical mass transfer among the constituents.

The starting point has been the kinetic model for reactive gases {proposed}
by Desvillettes, Monaco and Salvarani \cite{DMS-EJM-05},
which has been studied here under the standard diffusive scaling.

The form of the limiting equations \eqref{eq:MS}
shows the influence of
both the chemical reaction and the polyatomic structure of the mixture
on the evolution of the number densities of the constituents, as well as on the evolution
of the energy of the mixture.
Moreover, the evolution equation of the momentum of each constituent also
shows the influence of the polyatomic structure of the mixture, through the diffusion coefficients $D_{ij}$,
but it is not affected by the chemical reaction.
This is due to the fact that the considered chemical regime corresponds to
a slow reaction, in which the reactive process, in comparison with diffusion, has a small effect in the evolution.

\medskip

We highlight that the target equations obtained here
take into account the effect of chemical reactions, which impose to work in the non-isothermal setting.
Notice that the sum over the index $i$ of the right hand sides of the four (vectorial) equations appearing in the second line of the system \eqref{eq:MS} vanishes.
Therefore they do not constitute a set of four independent equations for the fluxes $\bfJ_i$.

The limit equations are compatible with the internal energy law of polyatomic gases and the law of mass action.
At the equilibrium, from both equations \eqref{eq:MS}
and the formulation of the source term \eqref{eq:Aend}, the law of mass action reads
$$
\frac{n_1 n_2}{n_3 n_4}
=\left(\frac{m_1 m_2}{m_3 m_4}\right)^{3/2}
                 \frac{q_1(T(t,\bfx)) q_2(T(t,\bfx))}{q_3(T(t,\bfx)) q_4 (T(t,\bfx))}
                 \exp \left( \frac{E}{k_B T(t,\bfx)} \right),
$$
which is exactly the same as in \cite{DMS-EJM-05}.

Of course, all the computations heavily depend on the choice of the functional forms of the weights $\varphi_i$
which can be -- in principle -- freely chosen, provided that they induce a macroscopic behavior consistent with the physical situation under investigation.

In case of polytropic gases, the weights $\varphi_i$ have the form $\varphi_i(I)=I^{\alpha_i}$,
where $\alpha_i=(k-2)\in\mathbb{N}$ and $k$ is the number of atoms of the species $i$. The special case $\alpha_i=0$ corresponds to diatomic molecules~\cite{DMS-EJM-05, Pavic-Ruggeri-Simic}.
Thanks to this choice of the {weights} $\varphi_i$ we can compute
$$
q_i(T(t,\bfx))=\int_{0}^{+\infty} I^{\alpha_i} e^{-I/(k_BT(t,\bfx))}\mathrm{d} I
= [k_B T(t,\bfx)]^{\alpha_i+1}(\alpha_i !)\,,
$$
with the convention $0!=1$, and analogously
$$
q_i^*(T(t,\bfx))=\int_{0}^{+\infty} I^{\alpha_i+1} e^{-I/(k_BT(t,\bfx))} \mathrm{d} I
= [k_B T(t,\bfx)]^{\alpha_i+2}(\alpha_i +1) !\,.
$$
The model is thus consistent with the energy law of polytropic molecules, which provides a linear dependence on temperature $T(t, \bfx)$.
More precisely, {the} total energy at the equilibrium state becomes
$$
\sum_{i=1}^4 n_i \left( \frac32 k_B T + \frac{q_i^*(T)}{q_i(T)} \right) = \sum_{i=1}^4 n_i \left( \frac32 + \alpha_i + 1 \right) k_B T(t,\bfx).
$$
Because of the total conservation of the mass during the chemical reaction (i.e. $\alpha_1+\alpha_2=\alpha_3+\alpha_4$), we can conclude that
$$
\frac{q_1(T(t,\bfx)) q_2(T(t,\bfx))}{q_3(T(t,\bfx)) q_4 (T(t,\bfx))}= \frac{(\alpha_1 !)(\alpha_2 !)}{(\alpha_3 !)(\alpha_4 !)} \, ,
$$
which is a constant and does not depend on the temperature. In particular, it is equal to one as soon as $(\alpha_1 !)(\alpha_2 !)=(\alpha_3 !)(\alpha_4 !)$.

Of course, more complicated {weights} $\varphi_i$ should be chosen to reproduce situations with non--polytropic gases, for which total energy is made up also by exponential functions of the global temperature, as it occurs for instance in the kinetic description involving a set of discrete internal energies for each gas~\cite{bis-rug-spi-18}.

\bigskip


\subsection*{Acknowledgments}

The paper is partially supported by the
Portuguese FCT Project UID/MAT/00013/2013,
by the PhD grant PD/BD/128188/2016,
by the bilateral Pessoa project
7854WM and 406/4/4/2017/S ``Derivation of macroscopic PDEs from kinetic theory (mesoscopic scale) and from interacting particle systems (microscopic scale)'', by the ANR project \textit{Kimega} (ANR-14-ACHN-0030-01), by the Italian Ministry of Education, University and Research (MIUR), \textsl{Dipartimenti di Eccellenza} Program - Dept. of Mathematics ``F. Casorati'', University of Pavia,
by the University of Parma and by the Italian National Institute of Higher Mathematics INdAM--GNFM.
\newpage

\section*{Appendix -- Proofs of the Properties}

In part {\it(A)} of this appendix, we give the central ideas to prove Proposition \ref{prop:BiRe} stated in Subsection~\ref{ssec:pre}.
Then, in parts {\it(B)} and {\it(C)}, we give some details of the computation of the production terms $A$ and $\bfB_i$ appearing in the limiting equations derived in Subsections \ref{ssec:mom0} and \ref{ssec:mom1}, respectively.

\medskip


\noindent
{\it (A) Proof of Proposition \ref{prop:BiRe}, Subsection \ref{ssec:pre}.}

\noindent
By expanding the quadratic term in the first exponential appearing in the Maxwelllian \eqref{eq:maxhip}, one can split it  into three exponentials.
Taylor expanding {with respect to $\varepsilon$} the resulting four exponentials as well as the other two {terms} in front
of the exponentials,
we can write $f^\varepsilon_i$ in the form
\begin{align*}
f^\varepsilon_i =
\frac{n^\varepsilon_i}{q_i(T)}
    &
    \left( \! \frac{m_i}{2\pi k_BT} \! \right )^{\!\!3/2}
    \exp \left( -\frac{m_iv^2}{2k_BT}-\frac{I}{k_B T} \right)
    \\
    & \times \left[
    1+ \varepsilon\frac{m_i \bfv \!\cdot\! {\bfu}_i^\varepsilon}{k_BT} +
    \varepsilon\frac{ m_i v^2 T_i^\varepsilon}{2k_B T^2} +
    \varepsilon\frac{I T_i^\varepsilon}{k_B T^2} -
    \varepsilon \frac{q_i^*(T)}{k_B T q_i(T)} \frac{T_i^\varepsilon}{T} -
    \frac32 \varepsilon \frac{T_i^\varepsilon}{T}
    +
    O(\varepsilon^2) \right ]
\end{align*}
where $q_i^*(T)$ has been defined in \eqref{eq:qistar}, and similar expansions hold for $f^\varepsilon_j$, $f^{\varepsilon'}_i$, $f^{\varepsilon'}_j$
as well as for $f^{\varepsilon'}_k$, $f^{\varepsilon'}_l$.

\medskip

\noindent
To prove part {\it(a)}, we perform the products ${f}^{\varepsilon'}_i{f}^{\varepsilon'}_j$, $f^\varepsilon_if^\varepsilon_j$
and use the conservation of energy for elastic collisions given in \eqref{eq:conservation_eng_elastic},
to obtain
\begin{align}
\nonumber
{f}^{\varepsilon'}_i{f}^{\varepsilon'}_j&-f^\varepsilon_if^\varepsilon_j \\
\nonumber
&= M_i^\varepsilon M_j^\varepsilon
      \Bigg\{ \varepsilon
      \bigg[ \frac{{\bf u}^\varepsilon_i\cdot m_i({\bfv}'_i- \bfv)}{k_B T}+\frac{{\bf u}_j^\varepsilon
      \cdot m_j( {\bfv}'_j-{\bfv}_j)}{k_B T}+\frac{ T_i^\varepsilon}{k_B T^2}\bigg( \frac{ m_i(v'_i)^2}{2} - \frac{ m_i v^2}{2}+I_i'-I\bigg)
      \\
&      \qquad -\frac{T_j^\varepsilon}{k_B T^2} \bigg( \frac{ m_j v_j^2}{2}-\frac{m_j (v'_j)^2}{2}+I_j-I_j' \bigg)
     \bigg]+O( \varepsilon^2) \Bigg\}.
\label{eq:bispecie_gain_and_loss_1}
\end{align}
Using the conservation of momentum  \eqref{eq:conservation_mom_elastic}
for the two first addends within the brackets
and the conservation of energy \eqref{eq:conservation_eng_elastic}
for the next two addends,
we obtain the result stated in equation~\eqref{eq:gain_minus_loss_bispecies}.
Therefore, part {\it(a)} of Proposition \ref{prop:BiRe} is proven.

\medskip

\noindent
Similarly, to prove part {\it(b)}, we perform the products ${f}^{\varepsilon'}_k{f}^{\varepsilon'}_l$, $f^\varepsilon_if^\varepsilon_j$
and we obtain
\begin{align*}
\left( \frac{m_i m_j}{m_k m_l} \right )^{\!3}
     &\!\! f^{\varepsilon'}_k\! f^{\varepsilon'}_l-f^\varepsilon_i\! f^\varepsilon_j
     \\
     &= \left (\frac{m_i m_j}{m_k m_l} \right )^{\!3}
     \frac{n_k}{q_k(T)} \frac{n_l}{q_l(T)}
     \frac{(m_k m_l)^{3/2}}{(2\pi k_BT)^3}
     \exp\!  \left(\! -\frac{m_k(v_k')^2 + m_l{(v_l')^2}}{2k_BT}\!-\! \frac{I'_k + I'_l}{k_BT}\!\right)
     \\
     &
     \quad \times \Bigg[
     1+\varepsilon\frac{m_k{\bfu}_k^\varepsilon\!\cdot\! {\bfv}_k'}{k_BT}
     + \varepsilon\frac{ m_k T_k^\varepsilon (v_k')^2}{2k_B T^2}
     + \varepsilon\frac{T_k^\varepsilon I_k'}{k_B T^2}
     - \varepsilon \frac{T_k^\varepsilon}{T} \frac{q_k^*(T)}{k_BT q_k(T)}
     - \varepsilon \frac{3T_k^\varepsilon}{2T}
     \\
     & \quad
     + \varepsilon\frac{m_l{\bfu}_l^\varepsilon\!\cdot\! {\bfv}'_l}{k_BT}
     + \varepsilon\frac{ m_l T_l^\varepsilon (v_l')^2}{2k_B T^2} +
     \varepsilon\frac{T_l^\varepsilon I_l'}{k_B T^2}
      - \varepsilon \frac{T_l^\varepsilon}{T} \frac{q_l^*(T)}{k_BT q_l(T)}
     - \varepsilon \frac{3T_l^\varepsilon}{2T}
     + O(\varepsilon^2) \bigg]
     \\
     &
     - \frac{n_i}{q_i(T)} \frac{n_j}{q_j(T)}
     \frac{(m_i m_j)^{\!3/2}}{(2\pi k_BT)^3}
     \exp\! \left(\! -\frac{m_i(v)^2 + m_j(v_j)^2}{2k_BT}\!-\! \frac{I + I_j}{k_BT}\! \right)
     \\
     & \quad \times
     \Bigg[1+\varepsilon\frac{m_i{\bfu}_i^\varepsilon\!\cdot\! \bfv}{k_BT}
     + \varepsilon\frac{ m_i T_i^\varepsilon v^2}{2k_B T^2}
     + \varepsilon\frac{T_i^\varepsilon I}{k_B T^2}
     - \varepsilon \frac{T_i^\varepsilon}{T} \frac{q_i^*(T)}{k_BT q_i(T)}
     - \varepsilon \frac{3T_i^\varepsilon}{2T}
     \\
     &  \quad
     + \varepsilon\frac{m_j{\bfu}_j^\varepsilon\!\cdot\! {\bfv}_j}{k_BT}
     + \varepsilon\frac{ m_j T_j^\varepsilon v_j^2}{2k_B T^2}
     + \varepsilon\frac{T_j^\varepsilon I_j}{k_B T^2}
     - \varepsilon \frac{T_j^\varepsilon}{T} \frac{q_j^*(T)}{k_BT q_j(T)}
     - \varepsilon \frac{3T_j^\varepsilon}{2T}
     +O(\varepsilon^2)
     \bigg].
\end{align*}
Putting together terms of the same order in $\varepsilon$,
we obtain the result stated in equation \eqref{eq:postf-pref_reactive}.
Therefore, part {\it(b)} of Proposition \ref{prop:BiRe} is also proven.
\QEDB

\bigskip


\noindent
{\it {(B)} On the computation of the production term $A^\varepsilon$ appearing in Subsections \ref{ssec:mom0}.}

\bigskip

\noindent
The production rate corresponding to species $i=3$ is provided by
\begin{align}
A^\varepsilon
        & =  \int_{\mathbb{R}^3}\!\int_0^{+\infty}\!Q_3^{react}{\varphi_3(I)} \dd I \dd{\bfv}
        \nonumber
        \\[1mm]
        & = \int_{\mathbb{R}^3} \! \int_{\mathbb{R}^3}\! \int_0^{+\infty} \!
               \int_0^{+\infty} \!\! \int_0^1 \! \int_0^1 \! \int_{\mathbb{S}^2} \!
                  \left[ \left(\frac{m_3 m_4}{m_1 m_2} \right)^{\!\!3} \!
                  M_1^{\varepsilon'} M_2^{\varepsilon '} - M_3^\varepsilon M_4^\varepsilon\right ] \!
                  V^{\gamma-1}H_3
                  \label{eq:tildeA_epsilon}
                  \\[1mm]
           & \quad \times  \Phi^{react}(I,I_4,R,r)\,
                  \frac{2\cos\theta(1-R)}{\varphi_3(I) (m_3m_4)^{2} }
                  \; \dd \sigma \dd r \dd R \dd I_4 \dd I  \dd{\bfv}_4 \dd{\bfv}
                  \!+\! O(\varepsilon).
                 \nonumber
\end{align}
Using the form \eqref{eq:M} for the distribution function,
together with the conservation of total energy \eqref{eq:conservation_eng_chemical}
for reactive collisions, we obtain,
\begin{equation}
\begin{aligned}
 &\left( \! \frac{m_3 m_4}{m_1 m_2} \! \right)^{\!\!3} \!\!
     M_1^{\varepsilon'} \! M_2^{\varepsilon '} \!-\! M_3^\varepsilon M_4^\varepsilon
     \\[1mm]
     &= \frac{(m_3m_4)^{3/2}}{ \big(2\pi k_BT\big)^{3}} \left[\!
          \bigg( \frac{m_3m_4}{m_1m_2} \bigg)^{\!\!3/2} \!\!
          \frac{n_1^\varepsilon n_2^\varepsilon}{q_1(T) q_2(T)}
          \exp \left( -\, \frac{E}{k_B T} \! \right)
          - \frac{n_3^\varepsilon n_4^\varepsilon}
                {q_3(T) q_4(T)} \right]
                 \\[1mm]
     &  \qquad \times \exp \! \left( \! - \frac{m_3 v^2 + m_4(v_4)^2}{2k_B T} -
         \frac{I \!+\! I_4}{k_B T} \right) .
\label{eq:tildeA_epsilon1}
\end{aligned}
\end{equation}
Substituting expression \eqref{eq:tildeA_epsilon1} into Equation \eqref{eq:tildeA_epsilon},
then transforming to the center of mass velocity $\Xx = (m_3 \bfv + m_4 \bfv_4)/M$ (with $M = m_3+m_4$) and relative velocity $\Vv = \bfv_4 - \bfv$, and taking into account the trivial definition of the Heaviside function $H_3$,
we get
\begin{eqnarray}
A^\varepsilon \!\!&\!\!=\!\!&\!\!
                 \frac{(m_3m_4)^{3/2}}{ \big(2\pi k_BT\big)^{3}}
                 \left[\!
                 \bigg( \frac{m_3m_4}{m_1m_2} \bigg)^{\!\!3/2} \!\!
                 \frac{n_1^\varepsilon n_2^\varepsilon}{q_1(T) q_2(T)}
                 \exp \left( - \frac{E}{k_B T} \! \right)
                 - \frac{n_3^\varepsilon n_4^\varepsilon}
                {q_3(T) q_4(T)} \right]
                 \nonumber
                 \\[1mm]
        & &   \times \! \int_{\mathbb{R}^3} \int_{\mathbb{R}^3} \! \int_0^{+\infty} \!\! \int_0^{+\infty} \!\!
                 \int_0^1 \! \int_0^1 \! \int_{\mathbb{S}^2} \!
                 \exp \! \left( \! - \frac{MX^2+\mu_{34}V^2}{2k_B T}\right)
                 \exp \! \left( -\frac{I\!+\!I_4}{k_B T} \!\right)
                 \label{eq:tildeAetau}
                 \\[1mm]
        & &   \times \, V^{\gamma -1} {\Phi^{react}(I , I_4 , R , r)}
                 2 b^{react}(\cos\theta) (m_3 m_4)^{-2} {\cos \theta} \, (1-R) \,
                 \dd \sigma \dd r \dd R \dd I_4 \dd I d{\bf{X}} d{\bf{V}}.
                 \nonumber
\end{eqnarray}

\medskip

\noindent
Now, the integral in $\bf X$ results in
\begin{equation}
\int _{{\mathbb{R}}^3} \exp \left( - \frac{M X^2}{2k_BT} \right ) \, d{\bf{X}}
    =  \left( \frac{2\pi k_B T}{M} \right)^{3/2},
\label{eq:tildeB_bar2}
\end{equation}
since
\begin{align*}
\int _{{\mathbb{R}}^3} \exp \left( - \frac{M X^2}{2k_BT} \right) \, d{\bf{X}}
      &= \int_0^{{+\infty}}\int_0^{\pi} \int_0^{2\pi}X^2 \exp \left( - \frac{M X^2}{2k_BT} \right)
            \sin\theta \,d \phi\, d \theta\,d{{X}} \\
       &= 4\pi\int_0^{{+\infty}}X^2 \exp \left( - \frac{M X^2}{2k_BT} \right) \, d{{X}}.
\end{align*}
Evaluating the previous integral using definition \eqref{eq:int_gamma_function}
gives the result in \eqref{eq:tildeB_bar2}.

\bigskip

\noindent
Also, the integral in $\bf{V}$ results in
\begin{equation}
\int_{\mathbb{R}^3} V^{\gamma -1} \exp \! \left( \! - \frac{\mu_{34}V^2}{2k_B T}\right) d{\bf{V}}
      = 2\pi \left( \frac{2k_BT}{\mu_{34}} \right)^{(\gamma+2)/2}
         \Gamma \left( \frac{\gamma + 2}{2} \right).
         \label{eq:tildeintV}
\end{equation}
Indeed,
\begin{align*}
\int_{\mathbb{R}^3} V^{\gamma -1} \exp \! \left( \! - \frac{\mu_{34}V^2}{2k_B T}\right)
      d{\bf{V}}
      &= \int_0^{{+\infty}}\int_0^{\pi} \int_0^{2\pi}V^{\gamma +1}
      \exp \! \left( \! - \frac{\mu_{34}V^2}{2k_B T}\right)
            \sin\theta \,d \phi\, d \theta\,dV\\
      &= 4 \pi\int_0^{{+\infty}} V^{\gamma +1}
      \exp \! \left( \! - \frac{\mu_{34}V^2}{2k_B T}\right) dV.
\end{align*}
Evaluating the previous integral using definition \eqref{eq:int_gamma_function}
gives the result in \eqref{eq:tildeintV}.

\medskip

\noindent
Finally, substituting equations \eqref{eq:tildeB_bar2} and \eqref{eq:tildeintV} into Equation \eqref{eq:tildeAetau},
we obtain the desired expression for $A^\varepsilon$ given in \eqref{eq:Aetauu}.

\medskip

\noindent
{\it {(C)} On the computation of the simple form of the production term $\bfB_i$
           appearing in Equation \eqref{eq:Bbar_final_1}} of Subsection \ref{ssec:mom1}.

\medskip

\noindent
To prove part {\it {(C1)}} it is enough to show that
\begin{equation}
\int_{\mathbb{S}^2}  \sigma {b_{ij}}(\cos\theta)\cos\theta \dd \sigma=0.
\label{eq:E_1}
\end{equation}
To see this, observe that
\begin{align*}
\int_{\mathbb{S}^2} & \sigma {b_{ij}}(\cos\theta)\cos\theta \dd \sigma \\
     &=  \int_0^\pi\int_0^{2\pi} \big( {{\hat{\bfx}}}
     \sin\theta \cos\phi+{\hat{\bfy}}\sin\theta \sin\phi +
     {\hat{\bfz}}\cos\theta \big) {b_{ij}}(\cos\theta)\cos\theta \sin\theta
     \; d\phi d\theta
     \\
     &=  {{\hat{\bfx}}} \! \int_0^\pi \! {b_{ij}}(\cos\theta)
     \sin^2\theta \cos\theta d\theta\underbrace{\int_0^{2\pi}\cos\phi d\phi}_{0}+{{\hat{\bfy}}}\int_0^\pi
     {b_{ij}}(\cos\theta) \sin^2\theta \cos\theta d\theta \! \underbrace{\int_0^{2\pi} \! \sin\phi d\phi}_{0}
     \\
     &+ {{\hat{\bfz}}} \!
     \underbrace{\int_0^\pi {b_{ij}}(\cos\theta) \sin\theta \cos^2\theta d\theta}_{0, \,\,\,{\rm since}\,\,
     {b_{ij}}(\cos\theta)\,\, \text{\rm is odd}}\int_0^{2\pi} \! d\phi.
\end{align*}

\medskip


\noindent
To prove part {\it {(C2)}} we first {transform to the centre of mass velocity and relative velocity  and obtain}
\begin{align}
\int_{\mathbb{R}^3} \!
      & \! \int_{\mathbb{R}^3}
      (\ba_{ij}\cdot {\bf{V}}) \bfv \, M_i^\varepsilon M_j^\varepsilon V^{\gamma-1} \dd{\bfv}_j \dd\bfv
      \nonumber \\
      &\!\!{=\!\frac{n_i^\varepsilon n_j^\varepsilon}{q_iq_j}\frac{{ (m_im_j)^{3/2}}}{(2\pi k_BT)^3}
      \exp \! \left( \! -\frac{I \!+\!I_j}{k_B T} \! \right ) \!\!
       \! \int_{\mathbb{R}^3}\!\!\int_{\mathbb{R}^3} \!\!({\ba}_{ij} \!\cdot\! {\bf{V}})\!
      \Bigg(\!{\bf{X}}\!-\!\frac{m_j{\bf{V}}}{m_i+m_j}\! \Bigg)\!\exp \! \left(\! -\frac{MX^2\!+\!\mu_{ij}{{V}}^2}{2k_BT} \! \right)\!V^{\gamma-1}\! d{\Vv} d{\bf{X}}}
      \nonumber\\
      & = \! \frac{n_i^\varepsilon n_j^\varepsilon}{q_iq_j}\frac{ (m_im_j)^{3/2}}{(2\pi k_BT)^3}
      \exp \! \left( \! -\frac{I \!+\!I_j}{k_B T} \! \right ) \!\!
      \Bigg[ \! \int_{\mathbb{R}^3} \!\! {\bf{X}} \! \exp \! \left( \! -\frac{M{{X}}^2}{2k_BT} \! \right) \! d{\bf{X}} \!
      \int_{\mathbb{R}^3} \!\!({\ba}_{ij} \!\cdot\! {\bf{V}})\!
      \exp \! \left(\! -\frac{\mu_{ij}{{V}}^2}{2k_BT} \! \right)\!V^{\gamma-1}\! d{\Vv}
      \nonumber \\
      & - \! \int_{\mathbb{R}^3} \! \exp \! \left(\! -\frac{M{{X}}^2}{2k_BT} \right) \! d{\bf{X}}\,
      \frac{m_j}{m_i+m_j}
      \int_{\mathbb{R}^3} \! {\bf{V}}({\ba}_{ij}\cdot {\bf{V}})
      \exp \! \left( - \! \frac{\mu_{ij}{{V}}^2}{2k_BT} \right)V^{\gamma-1} d{\Vv} \Bigg].
      \label{eq:int_B_bar}
\end{align}
Next we compute the integrals in $\bf X$ and $\bf V$ on the right hand side of {Eq.~\eqref{eq:int_B_bar}.}
We use the fact that any vector can be written in terms of a unit vector,
{\it i.e.} ${\bf{V}}=V\vec{\bfv}$, ${\bf{X}}=X\vec{\bfx}$,
where $X\!=\!\left \| {\bf{X}} \right \|$, $V\!=\!\left \| {\bf{V}} \right \|$.
Also, we expand the unit vectors $\vec{\bfv}$ and $\vec{\bfx}$
in terms of the Cartesian unit vectors $\hat\bfx$, $\hat\bfy$, $\hat\bfz$ in $\RR^3$.
Finally, we use the notation
${\ba}_{ij} \cdot  {\bf{V}} = \left\| {\ba}_{ij} \right \| \left \| {\bf{V}}  \right\| \cos\theta = a_{ij} V \cos\theta$.
\begin{enumerate}
\item[(a)] For the first integral in $\bf X$ appearing {in the last equality of} \eqref{eq:int_B_bar}, we obtain
\begin{equation}
\int_{\mathbb{R}^3} {\bf{X}} \exp \left(-\frac{M{{X}}^2}{2k_BT} \right )d{\bf{X}} =  0,
\hspace*{5cm}
\label{eq:B_bar1}
\end{equation}
since
\begin{align*}
\int_{\mathbb{R}^3} {\bf{X}}
    & \exp \left(-\frac{M{{X}}^2}{2k_BT} \right )d{\bf{X}}
    \!=\! \int_{\mathbb{R}^3} {{X}}\vec{{\bfx}}\, \exp\left(-\frac{M{{X}}^2}{2k_BT} \right ) d{\bf{X}}
    \\[1mm]
    & \!=\! \int_0^{{+\infty}} \!\! \int_0^\pi \!\! \int_0^{2\pi} \!\!
    X \big( {\hat{\bfx}} \sin\theta \cos\phi + {\hat{\bfy}} \sin\theta \sin\phi +
    {\hat{\bfz}} \cos\theta \big) \! \exp \! \left(\!-\frac{M{{X}}^2}{2k_BT} \!\right )\!X^2
    \sin\theta \, d\phi d\theta d{{X}}
    \\[2mm]
    & \!= {\hat{\bfx}}\int_0^{{+\infty}} \!\! X^{3} \exp \left( - \frac{M{{X}}^2}{2k_BT} \right ) dX
    \int_0^\pi \sin^2\theta \, d\theta \underbrace{\int_0^{2\pi}cos\phi \, d\phi}_{0}
    \\
    & \qquad + {\hat{\bfy}}\int_0^{{+\infty}} \!\!  X^{3} \exp\left(-\frac{M{{X}}^2}{2k_BT} \right ) dX
    \int_0^\pi \sin^2\theta d\theta \underbrace{\int_0^{2\pi} \sin\phi d\phi}_{0}
    \nonumber \\
    & \qquad \qquad + {\hat{\bfz}} \int_0^{{+\infty}} \!\!  X^{3} \exp \left( - \frac{M{{X}}^2}{2k_BT} \right ) dX
    \underbrace{\int_0^\pi \sin\theta \cos\theta d\theta}_{0}\int_0^{2\pi}d\phi
    = 0.
\end{align*}


\item[(b)] The second integral in $\bf X$ appearing {in the last equality of} \eqref{eq:int_B_bar} has been evaluated in part {{\it(B)}} of this Appendix, see Equation \eqref{eq:tildeB_bar2}.


\item[(c)]
For the first integral in ${\bf{V}}$ appearing {in the last equality of} \eqref{eq:int_B_bar}, we obtain
\begin{equation}
     \int_{\mathbb{R}^3} ({\ba}_{ij} \cdot {\bf{V}})
      \exp  \left(\! -\frac{\mu_{ij}{{V}}^2}{2k_BT}  \right)\!V^{\gamma-1} d{\Vv}=0
\label{eq:B_bar4}
\end{equation}
since
\begin{align*}
\int_{\mathbb{R}^3} \!\! ({\ba}_{ij} \!\cdot\! {\bf{V}})\exp \! \left( \! -\frac{\mu_{ij}{{V}}^2}{2k_BT} \! \right)V^{\gamma-1} d{\Vv}
      &\!=\! \int_{\mathbb{R}^3} a_{ij}V \cos\theta
      \exp  \left( -\frac{\mu_{ij}{{V}}^2}{2k_BT} \right)V^{\gamma-1} d{\Vv}
      \\
      &\!\!=\!a_{ij} \!\! \int_0^{{+\infty}}\!\!\!\int_0^{\pi}\!\!\!\int_0^{2\pi}\!\!\! V^{\gamma+2}\exp \!
      \left(\! -\frac{\mu_{ij}{{V}}^2}{2k_BT} \!\right) \sin\theta \cos\theta d\phi d\theta \!dV
      \\
      &=\!a_{ij} \!\! \int_0^{{+\infty}}\!\!\!V^{\gamma+2} \! \exp \! \left(\! -\frac{\mu_{ij}{{V}}^2}{2k_BT} \!\right) \! dV \!\!
      \underbrace{\int_0^{\pi}\!\!\! \sin\theta \cos\theta d\theta}_{0} \underbrace{\int_0^{2\pi}\!\!\! d\phi}_{2\pi}.
\end{align*}
Observe that the integral with respect to $V$ in the previous line can be evaluated using the expression \eqref{eq:int_gamma_function}.


\item[(d)] Lastly, for the second integral in $\bf V$ appearing {in the last equality of} \eqref{eq:int_B_bar}, we obtain

\begin{equation}
\int_{\mathbb{R}^3}({\ba}_{ij} \! \cdot \! {\bf{V}}){\bf{V}} \exp \! \left( \! - \frac{\mu_{ij}V^2}{2k_B T} \right) \! V^{\gamma-1} d{\bf{V}}
    =  \frac{2\pi}{3} {{\ba}_{ij}} \, \Gamma \! \left( \! \frac{\gamma+4}{2} \right) \!
        \left( \! \frac{2k_B T}{\mu_{ij}} \right )^{\!\! \frac{\gamma+4}{2}},
\label{eq:B_bar3}
\end{equation}
since
\begin{align*}
\int_{\mathbb{R}^3}({\ba}_{ij} & \! \cdot \! {\bf{V}}){\bf{V}}
       \exp \left( -\frac{\mu_{ij}V^2}{2k_B T} \right) V^{\gamma-1} d{\bf{V}} \\
       & = \int_{\mathbb{R}^3} (a_{ij} \cos\theta)V^2 V^{\gamma-1}  \vec{{\bfv}}
       \exp\left ( -\frac{\mu_{ij}V^2}{2k_B T} \right ) d{\bf{V}} \\
       & = \int_0^{{+\infty}} \!\! \int_0^{\pi} \!\! \int_0^{2\pi} \!\! (a_{ij} \cos\theta)V^2 \Big( {{\hat{\bfx}}}
       \sin\theta \cos\phi + {\hat{\bfy}} \sin\theta \sin\phi + {\hat{\bfz}} \cos\theta \Big)
       \hspace*{2cm}
       \\
       & \qquad\qquad \times V^{\gamma-1} \exp \left( \! -\frac{\mu_{ij}V^2}{2k_B T} \! \right) V^2 \sin\theta \, d\phi d\theta d{{V}} \\
       &=a_{ij} \, {\hat{\bfx}} \!\int_0^{{+\infty}}\!\! V^{\gamma+3}\!
       \exp\!\left (\! -\frac{\mu_{ij}V^2}{2k_B T} \!\right ) dV
       \int_0^{\pi}\!\! \sin^2\theta \cos\theta d\theta \underbrace{\int_0^{2\pi}\! \cos\phi \, d\phi}_{0} \\
       & \qquad\qquad + a_{ij}{\hat{\bfy}}\!\int_0^{{+\infty}} \!\! V^{\gamma+3}\!
       \exp\!\left (\! -\frac{\mu_{ij}V^2}{2k_B T} \!\right ) dV
       \int_0^{\pi}\!\! \sin^2\theta \cos\theta d\theta \underbrace{\int_0^{2\pi}\!\! \sin\phi \, d\phi}_{0} \\
       & \qquad\qquad + a_{ij}{\hat{\bfz}}\int_0^{{+\infty}}  V^{\gamma+3}
       \exp \left( - \frac{\mu_{ij}V^2}{2k_B T} \right) dV
       \underbrace{\int_0^{\pi} \sin\theta \cos^2\theta d\theta}_{\frac{2}{3}}\,\underbrace{\int_0^{2\pi} \, d\phi}_{2\pi} \\[-2mm]
       & = \frac{2\pi}{3} {{\ba}_{ij}} \Gamma \left( \frac{\gamma+4}{2} \right)
       \left( \frac{2k_B T}{\mu_{ij}} \right )^{\!\!\frac{\gamma+4}{2}} ,
\end{align*}
where we have used definition \eqref{eq:int_gamma_function} {for the integral with respect to $V$ } to obtain the last equality.
\end{enumerate}

\medskip

\noindent
To conclude the proof of {part {\it (C2)}}, it is enough now to substitute
expressions \eqref{eq:tildeB_bar2}, \eqref{eq:B_bar1}, \eqref{eq:B_bar4} and \eqref{eq:B_bar3}
into {equation} \eqref{eq:int_B_bar} and we obtain
expression  \eqref{eq:int_vivj_Bbar}.

\medskip


\noindent
To prove part {\it {(C3)}}, observe that using the definition of ${\bfv}_i'$ given in \eqref{eq:v_iv_jprime},
we obtain that
\begin{align*}
\frac{m_i(v'_i)^2}{2} & - \frac{ m_i v^2}{2} +I_i' -I
        \!=\! \bigg(\! \frac{m_i^3}{2(m_i\!+\!m_j)^2} \!-\!
        \frac{m_i}{2}\!\bigg) v^2 + \frac{m_j \mu_{ij}v_j^2}{2(m_i\!+\!m_j)} +
        \frac{m_i \mu_{ij} ({\bfv}\!\cdot\! {\bfv}_j)}{m_i+m_j}
        \\
        &+ \frac{m_i \mu_{ij} ({\bfv}\!\cdot\! \sigma)}{m_i+m_j}\!\bigg(\!\frac{2R\cal E}{\mu_{ij}}\! \bigg)^{\!\frac{1}{2}}
        +\frac{m_j \mu_{ij} ({\bfv}_j\cdot \sigma)}{m_i+m_j}\bigg(\frac{2R\cal E}{\mu_{ij}} \bigg)^{\frac{1}{2}}
        +\frac{m_j \mu_{ij}}{2(m_i+m_j)}\bigg(\frac{2R\cal E}{\mu_{ij}} \bigg)+I_i'-I.
\end{align*}
Substituting the previous expression into the third term on the right hand side of
Equation \eqref{eq:B_i*}
and then expanding, we obtain that such term can be given by
\begin{eqnarray}
\lefteqn{\sum_{\substack{j=1\\ j\neq i}}^{4}\frac{ T_i^\varepsilon-T_j^\varepsilon}{k_B T^2}\bigg(\! \frac{m_i^3}{2(m_i\!+\!m_j)^2}\!-\!\frac{m_i}{2}\!\bigg)
{\int _{{\mathbb{R}}^3}\!\int_{\mathbb{R}^3} \!\int_0^{+\infty}\!\!\! \int_0^{+\infty} \!\! \int_0^1 \! \int_0^1 \! \int_{\mathbb{S}^2}}  \!
    \bfv v^2 M_i^\varepsilon M_j^\varepsilon V^{\gamma-1} \Phi_{ij}(I,I_j,R,r)}
    \nonumber
    \\[-3mm]
& & \hspace*{6cm} \dps
    \times  2 \, b_{ij} (\cos\theta) \cos\theta(1-R) \,\dd \sigma \dd r \dd R \dd I_j \dd I \dd{\bfv}_j \dd\bfv
    \nonumber\\
& &  \dps +\sum_{\substack{j=1\\ j\neq i}}^{4}\frac{ T_i^\varepsilon-T_j^\varepsilon}{k_B T^2}\frac{m_j \mu_{ij}}{2(m_i\!+\!m_j)}{\int _{{\mathbb{R}}^3}\!\int_{\mathbb{R}^3} \!\int_0^{+\infty}\!\!\! \int_0^{+\infty} \!\! \int_0^1 \! \int_0^1 \! \int_{\mathbb{S}^2}} \!\bfv v_j^2  M_i^\varepsilon M_j^\varepsilon V^{\gamma-1} \Phi_{ij}(I,I_j,R,r)
    \nonumber\\[-3mm]
& &  \hspace*{6cm} \dps \times  2 \, b_{ij} (\cos\theta) \cos\theta(1-R) \,\dd \sigma \dd r \dd R \dd I_j \dd I \dd{\bfv}_j \dd\bfv
    \nonumber\\
& &  \dps +\sum_{\substack{j=1\\ j\neq i}}^{4}\frac{ T_i^\varepsilon-T_j^\varepsilon}{k_B T^2}\frac{m_i \mu_{ij}}{m_i+m_j}{\int _{{\mathbb{R}}^3}\!\int_{\mathbb{R}^3} \!\int_0^{+\infty}\!\!\! \int_0^{+\infty} \!\! \int_0^1 \! \int_0^1 \! \int_{\mathbb{S}^2}}  \bfv ({\bfv}\!\cdot\! {\bfv}_j)  M_i^\varepsilon M_j^\varepsilon V^{\gamma-1} \Phi_{ij}(I,I_j,R,r)
    \nonumber
    \\[-3mm]
& & \hspace*{6cm} \dps \times  2 \, b_{ij} (\cos\theta) \cos\theta(1-R) \,\dd \sigma \dd r \dd R \dd I_j \dd I \dd{\bfv}_j \dd\bfv
    \nonumber\\
& &  \dps +\sum_{\substack{j=1\\ j\neq i}}^{4}\frac{ T_i^\varepsilon\!-\!T_j^\varepsilon}{k_B T^2}\frac{m_i \mu_{ij} }{m_i+m_j}\!
    \bigg(\!\frac{2R\cal E}{\mu_{ij}}\! \bigg)^{\!\frac{1}{2}}{\!\!\int _{{\mathbb{R}}^3}\!\int_{\mathbb{R}^3} \!\int_0^{+\infty}\!\!\! \int_0^{+\infty} \!\! \int_0^1 \! \int_0^1 \! \int_{\mathbb{S}^2}} \! \bfv({\bfv}\!\cdot\! \sigma) M_i^\varepsilon M_j^\varepsilon V^{\gamma-1} \Phi_{ij}(I,I_j,R,r)
    \nonumber
    \\[-3mm]
& & \hspace*{6cm} \dps \times  2 \, b_{ij} (\cos\theta) \cos\theta(1-R) \,\dd \sigma \dd r \dd R \dd I_j \dd I \dd{\bfv}_j \dd\bfv
    \label{eq:E_3} \\
& &  \dps +\sum_{\substack{j=1\\ j\neq i}}^{4}\frac{ T_i^\varepsilon\!-\!T_j^\varepsilon}{k_B T^2}\frac{m_j \mu_{ij} }{m_i+m_j} \!
    \bigg(\!\frac{2R\cal E}{\mu_{ij}}\! \bigg)^{\!\frac{1}{2}}{\!\!\int _{{\mathbb{R}}^3}\!\int_{\mathbb{R}^3} \!\int_0^{+\infty}\!\!\! \int_0^{+\infty} \!\! \int_0^1 \! \int_0^1 \! \int_{\mathbb{S}^2}}  \bfv({\bfv}_j\!\cdot\! \sigma) M_i^\varepsilon M_j^\varepsilon V^{\gamma-1} \Phi_{ij}(I,I_j,R,r)
    \nonumber
    \\[-3mm]
& & \hspace*{6cm}\times  2 \, b_{ij} (\cos\theta) \cos\theta(1-R) \,\dd \sigma \dd r \dd R \dd I_j \dd I \dd{\bfv}_j \dd\bfv
    \nonumber\\
& &  \dps +\sum_{\substack{j=1\\ j\neq i}}^{4}\frac{ T_i^\varepsilon-T_j^\varepsilon}{k_B T^2}\frac{m_j \mu_{ij} }{2(m_i+m_j)}\!\bigg(\!\frac{2R\cal E}{\mu_{ij}}\! \bigg){\int _{{\mathbb{R}}^3}\!\int_{\mathbb{R}^3} \!\int_0^{+\infty}\!\!\! \int_0^{+\infty} \!\! \int_0^1 \! \int_0^1 \! \int_{\mathbb{S}^2}}  \bfv M_i^\varepsilon M_j^\varepsilon V^{\gamma-1} \Phi_{ij}(I,I_j,R,r)
    \nonumber
    \\[-3mm]
& & \hspace*{6cm} \dps \times  2 \, b_{ij} (\cos\theta) \cos\theta(1-R) \,\dd \sigma \dd r \dd R \dd I_j \dd I \dd{\bfv}_j \dd\bfv
    \nonumber\\
& &  \dps +\sum_{\substack{j=1\\ j\neq i}}^{4}\frac{ T_i^\varepsilon-T_j^\varepsilon}{k_B T^2}{\int _{{\mathbb{R}}^3}\!\int_{\mathbb{R}^3} \!\int_0^{+\infty}\!\!\! \int_0^{+\infty} \!\! \int_0^1 \! \int_0^1 \! \int_{\mathbb{S}^2}}  \bfv I_i' M_i^\varepsilon M_j^\varepsilon V^{\gamma-1} \Phi_{ij}(I,I_j,R,r)
    \nonumber\\[-3mm]
& & \hspace*{6cm} \dps \times  2 \, b_{ij} (\cos\theta) \cos\theta(1-R) \,\dd \sigma \dd r \dd R \dd I_j \dd I \dd{\bfv}_j \dd\bfv
    \nonumber\\
& &  \dps -\sum_{\substack{j=1\\ j\neq i}}^{4}\frac{ T_i^\varepsilon-T_j^\varepsilon}{k_B T^2}{\int _{{\mathbb{R}}^3}\!\int_{\mathbb{R}^3} \!\int_0^{+\infty}\!\!\! \int_0^{+\infty} \!\! \int_0^1 \! \int_0^1 \! \int_{\mathbb{S}^2}}  \bfv I M_i^\varepsilon M_j^\varepsilon V^{\gamma-1} \Phi_{ij}(I,I_j,R,r)
    \nonumber\\[-3mm]
& &  \hspace*{6cm} \dps
    \times  2 \, b_{ij} (\cos\theta) \cos\theta(1-R) \,\dd \sigma \dd r \dd R \dd I_j \dd I \dd{\bfv}_j \dd\bfv.
        \nonumber
\end{eqnarray}
Transforming the above integrals into the center of mass velocity and relative velocity and then evaluating the resulting integrals
using the same strategy as in part ${\it {(C2)}}$, we obtain that the six fold integrals in ${\bfv}$ and ${\bfv}_j$ vanish
in the first, second, third, sixth, seventh and eight terms of the above expression \eqref{eq:E_3}.
Furthermore,
using Equation \eqref{eq:E_1}, we conclude that the fourth and fifth terms of expression \eqref{eq:E_3} also vanish.


\newpage
$\;$



\end{document}